\documentclass{article}
\usepackage{graphicx} 

\title{JoLT2}
\date{Nov 28, 2026}
\usepackage[utf8]{inputenc}
\usepackage[a4paper, total={6in, 9in}]{geometry}
\usepackage{hyperref}
\usepackage{setspace}
\usepackage[all,cmtip]{xy}
\usepackage{hyperref}
\usepackage{tikz}
\usepackage{tikz-cd}
\usepackage{cancel}
\usepackage{ulem}
\usepackage[all]{xy}
\usepackage{tikz-cd}
\hypersetup{
    linktoc=all,     
    linkcolor=black,  
}

\usepackage{inputenc, mathtools, amssymb, graphicx,mathrsfs}
\usepackage{xcolor}
\DeclareMathOperator{\Spec}{Spec}

\newcommand{\C}{\mathbb{C}}

\newcommand{\PP}{\mathbb{P}}
\newcommand{\Z}{\mathbb{Z}}

\newcommand{\g}{\mathfrak g}

\newcommand{\Lie}{\operatorname{Lie}}

\newcommand{\eat}[1]{}
\newcommand{\lie}[1]
{\mathsf{#1}}
\newcommand{\wmf}[1]{\widehat{\mf{#1}}}
\newcommand{\wl}[1]{\widehat{\lie{#1}}}
 \usepackage{thm-restate}
 \usepackage{ytableau}
\usepackage{ulem}
\usepackage{authblk}
\title{Cancellation algebras and groups associated with orbit closures. }
\author[1]{Bharat Adsul}
\author[1]{Milind Sohoni}
\author[2]{K Venkata Subrahmanyam}
\affil[1]{Department of Computer Science, IIT Mumbai, Powai 400076 }
\affil[2]{Chennai Mathematical Institute, Siruseri, Chennai 600103}

\newtheorem{theorem}{Theorem}[section]
\newtheorem{lemma}[theorem]{Lemma}
\newtheorem{prop}[theorem]{Proposition}
\newtheorem{corollary}[theorem]{Corollary}
\newtheorem{remark}[theorem]{Remark}
\newtheorem{ex}[theorem]{Example}
\newtheorem{defn}[theorem]{Definition}

\newtheorem{assume}[theorem]{Assumption}
\newcommand{\ignore}[1]{}
\newcommand{\mf}[1]{\mathfrak{#1}}

\begin{document}
\maketitle
\begin{abstract}
Let $\gamma(t)y\to z$ be a one-parameter degeneration
in a finite-dimensional representation of an algebraic group $G$. Assume that $\gamma(t)$ is ample.
We construct a descending chain of cancellation Lie algebras,
from the ambient Lie algebra through the stabilizer of $z$
to the limiting stabilizer of $y$.
We show that there is a corresponding filtration of the limiting orbit tangent space
with successive quotients naturally isomorphic to quotients
of the cancellation algebras as modules.
Using jet stabilizers, we realize these algebras as Lie algebras of algebraic subgroups of the stabilizer of $z$.

Examples involving matrix actions and binary forms illustrate
the construction and its dependence on the path of degeneration.
\end{abstract}

\section{Introduction}
Let $X$ be a finite dimensional vector space over $\C$ and let $G \subseteq GL(X)$ be an algebraic group defined over $\C$. 
Let $V$ be a finite dimensional $G$-module via a rational map $\rho : G \rightarrow GL(V)$ defined over $\C$. For a $g\in G$ and $v\in V$, let $\rho (g)\cdot v$, or simply $gv$  denote the action of $g$ on $v$ via $\rho$. 

Let ${\cal O}=\C [[t]]$ be the ring of formal power series with coefficients in $\C$. Let $k = \mathbb{C}((t))$ denote its quotient field of formal Laurent series 
in $t$ with coefficients in $\mathbb{C}$, whose elements may have 
a pole of finite order at $t = 0$. 
For any algebraic subgroup $L\subseteq G$ defined over $\C$, let $L(k) $ be elements of $L$ with entries in $k$ and $L({\cal O})$ be those elements with entries in ${\cal O}$. 

Let $y\in V$ and $O(y)$ denote the $G$-orbit of $y$ and $\overline{O(y)}$, its Zariski closure in $V$. Let $z\in \overline{O(y)}$ but $z\not \in O(y)$. Moreover, let $\gamma (t) \in G(k)$ be a one-parameter family such that:

\begin{equation} \label{eqn1ps}  y(t)=\gamma(t)\cdot y = \sum_{i\geq 0} t^i y_i = z + t^e y_e +\ldots 
\end{equation}
where $z=y_0$ and $y_e \neq 0$.  We denote this as $y\stackrel{\gamma}{\rightarrow}z$ and call $z$ the leading term or limit of $y$ under $\gamma(t)$ and $y_e$ as the tangent of approach. We say $\gamma(t)$ is a one parameter family (1-PF) driving $y$ to $z$. It is known that  $z\in \overline{O(y)}$ implies the existence of such a $\gamma (t)\in G(k)$. When $G$ is reductive, $\gamma (t)$ assumes a special form, see \cite{iwahori1965some, mumford1994geometric}.

The motivation of this paper (and our earlier paper, \cite{ASS25jolt}) is to study the following question. Given elements $z, y \in V$ with large distinctive stabilizers $H,K\subseteq G$, and the condition $y\stackrel{\gamma}{\rightarrow}z$, utilize $\gamma (t)$ to connect the stabilizers $K=G_y$ and $H=G_z$. \eat{One formulation of the above would be to locate intermediate $G$-varieties lying between $\overline{O(z)}$ and $\overline{O(y)}$ which carry groups intermediate to $H$ and $K$.}

We begin with some notation. Let ${\mf{h}}$ and ${\mf{k}}$ be the Lie algebras of $H$ and $K$. 
Let $T_z O(z) ={\mf{g}}\cdot z$ (or simply $T_z$), denote the tangent space of the orbit $O(z)$ at $z$. Note that $T_z $ is an $H$-module. Let $\overline{N}$ be the $H$-module (as well as the ${\mf{h}}$-module) $V/T_z $ which corresponds to a ``normal'' slice at $z$ to the tangent space $T_z $. Let $\overline{y_e}\in \overline{N}$ denote the coset $y_e +T_z$.

One approach in \cite{ASS25jolt} was to use $\gamma$ to construct a {\bf leading term Lie algebra} $\widehat{\mf{k}(t)}$ (or simply $\widehat{\mf{k}}$) from $\mf{k}$, from the Lie algebra of $K$. It was shown that  $dim(\widehat{\mf{k}})=dim(\mf{k})$ and that $\widehat{\mf{k}} \subseteq \mf{h}_{\overline{y_e}}$, where $\mf{h}_{\overline{y_e}}$, the Lie algebra stabilizer of $\overline{y_e}\in \overline{N}$. 

This was further connected to the local geometry of $z\in \overline{O(y)}$ through the use of the graded algebra associated with the ideal of $O(z)$. 
Let $I_y \subseteq I_z \subseteq \C[V]$ be the ideals of $\overline{O(z)}$ and $\overline{O(y)}$ and let $R$ be the associated graded ring $R=\oplus_{i\geq 0} I_z^i/I_z^{i+1}$. 
The spectrum $Spec(R)$, which is also a $G$-variety, consists of an open set $W$ isomorphic to $G\times^H \overline{N}$.\eat{The model within $V$ of this is of course $W = \{ (z', \overline{v'}| z'\in O(z), v' \in V/T_{z'} \}$ (where $T_{z'}$ is the tangent space at $z'$ to $O(z)$.} 
$I_y$ gives us an ideal $\overline{I}_y\subseteq R$, the ring $R_y=R/\overline{I}_y$, and the set $W_y\subseteq Spec(R_y) \cap W$ which is open in $Spec(R_y)$. Then, there is a closed $H$-variety $Z_y \subset \overline{N}$ such that $W_y \cong G\times^H Z_y$. 
 Thus, the embedding $\overline{O(z)}\rightarrow \overline{O(y)}$ is locally described by the containment of the zero section $O(z)\rightarrow W_y$ and the reverse containment of Lie algebras $\mf{h}\supseteq \wmf{k}$. 
 
This paper obtains a refinement of the containment $\mf{h}\supseteq \wmf{k}$ by  constructing a chain of subalgebras obtained from the path $\gamma$.

\noindent
We now state our results. They assume the {\it ampleness} of $\gamma$, see \ref{def:ample}. 

\begin{defn}
    Let $y\stackrel{\gamma}{\rightarrow} z$ with $z=y_0$. Then:
    \begin{enumerate} 
    \item For an element $\lie{g}\in \mf{g}$ and $\lie{g}\neq 0$, let $\lie{g}(t)=\gamma(t) \lie{g} \gamma (t)^{-1} = t^a \lie{g}_a+ \ldots$ with $\lie{g}_a \neq 0$. We define $d(\lie{g})=a$ and $\widehat{\lie{g}}$ as $\lie{g}_a$. For $0\in \mf{g}$, $\widehat{0}=0$, while $d(0)$ is taken to be $\infty$. 
    \item For a $v\in V, v\neq 0$ and $v(t)=\gamma (t)\cdot v = \sum_{j\geq b} t^j v_j$, with $v_b \neq 0$, we define $d(v)$ as $b$ and $\widehat{v}$ as $v_b$. For $0\in V$, $\widehat{0}=0$, while $d(0)$ is taken to be $\infty$. 
    \end{enumerate}
\end{defn}

\begin{defn} \label{def:ample}
    Let $M$ be a finite dimensional $G$-module and $N\subseteq M$ be a subspace. Let $\widehat{N}$ be the vector space generated by $\{ \widehat{n}| n\in N\}$. We say $N$ is ample is $dim(\widehat{N})=dim(N)$. We say $\gamma (t)$ is ample, if for every pair $(M,N)$ as above, $N$ is ample. 
\end{defn}
In general, we have $dim(\widehat{N})\leq dim(N)$. All one parameter subgroups are ample. 

Note that for any non-zero element $\lie{g}\in \mf{g}$ and $v\in V$, $\gamma (t) \lie{g} v= \lie{g}(t) v(t)$, and hence $d(\lie{g}v)\geq d(
\lie{g})+d(v)$. Hence, for any non-zero $\lie{g} \in \mf{g}$, since $d(y)= 0$, we have $d(\lie{g}y)\geq d(\lie{g})$. 

\begin{defn} \label{defn:delta}
    For a non-zero Lie algebra element $\lie{g}\in \mf{g}$, we define its {\bf cancellation number} $\delta_y (\lie{g})$ (or simply $\delta (\lie{g})$) as $d(\lie{g}y)- d(\lie{g})$.
\end{defn}

\begin{theorem} \label{thm:filter}
    With the the above notation, we have:
    \begin{enumerate} \item Let $\wmf{g}_r$ be the vector space spanned by all elements $\{ \widehat{\lie{g}}|\delta(\lie{g})\geq r \}$. Then $\wmf{g}_r$ is a Lie subalgebra. 
    \item Moreover, there is an $s$ such that $\mf{g}=\wmf{g}_0 \supseteq \wmf{g}_1 \supseteq \ldots \supseteq \wmf{g}_s\supsetneq \wmf{g}_{s+1}= \widehat{\kappa} $,  with $ \wmf{g}_1 =\mf{h}$ and $\wmf{g}_2 = \mf{h}_{\overline{y_e}}$. 
     \item  Let $Q_r \subseteq V$ be the vector space generated by the set $\{ \wl{gy}|\delta(\lie{g})\leq r\}$. Then $Q_r \subseteq \widehat{T_y}$ and is a $\wmf{g}_{r+1}$-module. Moreover, $Q_0 \subseteq Q_1 \ldots \subseteq Q_s$ and $Q_r/Q_{r-1} \cong \wmf{g}_r / \wmf{g}_{r+1}$ as a $\wmf{g}_{r+1}$-module.
        \item Finally, $Q_0= T_z O(z)\cong \mf{g}/\mf{h}$ and $Q_s =\widehat{Ty}\cong \mf{g}/\widehat{\kappa}$.
        \end{enumerate}
\end{theorem}

We call $(\wmf{g}_r )$ the cancellation algebras for the limit $y\stackrel{\gamma}{\rightarrow}z$.

Next, we come to the corresponding groups. Let $\C^r[t]= \C[t]/(t^{r+1})$. Let $V'=V\otimes {\cal O}$ and $V^r=V\otimes \C^r [t]=V'/t^{r+1}V'$. We also use the notation $V^{\infty}$ for $V'$. Then we have the sequence:
\[ V^\infty \stackrel{\pi^s}{\rightarrow} V^s \ldots  \stackrel{\pi^{i+1}}{\rightarrow} V^{i+1}
\stackrel{\pi^i}{\rightarrow} V^i \ldots 
\stackrel{\pi^1}{\rightarrow} V^1
\stackrel{\pi^0}{\rightarrow} V^0 \cong V\]
Let $G'=G({\cal O})$ and $G^r=G(\C^r [t])$, then we have the natural surjections:
\[ G' \stackrel{\pi^s}{\rightarrow} G^s \ldots  \stackrel{\pi^{i+1}}{\rightarrow} G^{i+1}
\stackrel{\pi^i}{\rightarrow} G^i \ldots 
\stackrel{\pi^1}{\rightarrow} G^1
\stackrel{\pi^0}{\rightarrow} G^0 \cong G\]

Note that, via $\rho$, $G'$ acts on $V^r$ for all $r$ and this factorizes through $G^r$. Thus, thus the action of $G^r$ on $V[r]$ matches that of $\pi^r (G')$ on $V[r]$. Also note that for $s>r$, $\pi^r :V^s \rightarrow V^r$ is $\rho$-equivariant. 

\eat{In other words we have the following commutative diagram:
\[ \begin{array}{ccccc}
G' &\times& V^\infty & \stackrel{\rho}{\longrightarrow} & V^\infty \\
\downarrow \pi^s & & \downarrow \pi^s& &\downarrow \pi^s \\
G^s &\times& V^s & \stackrel{\rho}{\longrightarrow} & V^s \\
\downarrow \pi^r & & \downarrow \pi^r& &\downarrow \pi^r \\
G^r &\times& V^r & \stackrel{\rho}{\longrightarrow} & V^r \\
\end{array}
\]
}
\begin{theorem} \label{theorem:grouptower}
Consider the element $y(t)\in V^\infty$ and let $y^r=\pi^r (y(t))$. We have the equivariant sequence:
\[ y(t)=y^\infty \stackrel{\pi^s}{\rightarrow} y^s \ldots  \stackrel{\pi^{i+1}}{\rightarrow} y^{i+1}
\stackrel{\pi^i}{\rightarrow} y^i \ldots 
\stackrel{\pi^1}{\rightarrow} y^1
\stackrel{\pi^0}{\rightarrow} y^0 =y_0 \in  V\]
Let $H^{r+1}=G^r_{y^r}$, the stabilizer of $y^r$ in $G^r$. 
\begin{enumerate}
    \item Then we have the sequence of homomorphisms: 
\[ H^\infty \stackrel{\pi^s}{\rightarrow} H^{s+1} \ldots  \stackrel{\pi^{i+1}}{\rightarrow} H^{i+1}
\stackrel{\pi^i}{\rightarrow} H^{i+1} \ldots 
\stackrel{\pi^1}{\rightarrow} H^2
\stackrel{\pi^0}{\rightarrow} H^1 =H =G_z  \]

\item Let $H^r (0)=\pi^0 (H^r)$ be the evaluation of $H^r$ at $t=0$. 
These are algebraic
subgroups and form a descending chain
\[
H=H^1(0)\supseteq H^2(0)\supseteq\cdots.
\]
Moreover, $\Lie\bigl(H^r(0)\bigr)=\wmf{g}_r$.
In particular, the algebraic subgroup $\widehat K:=H^{s+1}(0)\subseteq H$
satisfies 
$\Lie(\widehat K)=\wmf{k}.$
For every $r\geq s+1$, we have
\[\Lie\bigl(H^r(0)\bigr)=\wmf{k},
\qquad
\bigl(H^r(0)\bigr)^\circ=\widehat K^\circ.
\]
\item Moreover, the kernel $\pi_r^* : T_{y^{r+1}} O(y^{r+1}) \rightarrow T_{y^r} O(y^r)$ is isomorphic to the space $Q_r$. 
\eat{\item The map $-d:\mf{g}\rightarrow \Z$ is convex. Moreover,  there is bounded subgroup $K(t)\subseteq G(k)$ (in the Bruhat-Tits sense) corresponding to the data $(\mf{k}, -d)$ such that $\gamma (t) K(t) \gamma^{-1} (t)=(H')^s\subseteq G'$. }
\end{enumerate}
    
\end{theorem}

\section{Cancellation Algebras}
For Lie algebras we refer to \cite{humphreys2012introduction}, for algebraic groups we refer to \cite{humphreys2012linear}. We need some elementary notions from the theory of jets for which we use \cite{ein2008jetschemessingularities}. 

Recall that $\gamma (t) y = y_0 t^0 + y_e t^e +\ldots $, with $y_e \neq 0$. This is denoted by $y\stackrel{\gamma}{\rightarrow} y_0$ or $y\stackrel{\gamma}{\rightarrow} (y_0, y_e)$. In this paper, we will assume that $e=1$, and $y_0 =z$ as above. We also assume that $\gamma$ is ample.

Ampleness is satisfied when $\gamma (t)$ is a 1-parameter subgroup. Otherwise, it appears to be a delicate question, even when $G$ is a reductive group. 

\subsection{The bi-filtration}
\begin{defn}
    \label{defn:bifil}
    Let $\mf{g}^i =\{ \lie{g}\in \mf{g}| d(g)\geq i\}$ and $\mf{g}^{(j)}= \{ \lie{g}\in \mf{g}| d(\lie{g}y)\geq j \}$. Note that $\mf{g}^i$ and $\mf{g}^{(j)}$ are vector subspaces of $\mf{g}$. We define $\mf{g}^{i,j}$ as $\mf{g}^i \cap \mf{g}^{(j)}$.
\end{defn}

\begin{lemma}
    \label{lemma:bifil}
    There is an $a$ such that $\mf{g}^a= \mf{g}$ and that for all $i$, $\mf{g}^i  \supseteq \mf{g}^{i+1}$. Similarly, there is a $b\geq a$ such that $\mf{g}^{(b)}= \mf{g}$ and that for all $j$, $\mf{g}^{(j)}  \supseteq \mf{g}^{(j+1)}$. There is a smallest $s$ such that $\mf{g}^{(r)}=\mf{k}$ for all $r> s$. 
\end{lemma}

\noindent
{\bf Proof}: Only the last assertion needs proof.  We know that $\mf{g}^{(i)}$ are subspaces of $\mf{g}$. Since $\lie{k}y=0$ for $\lie{k}\in \mf{k}$, we must have $\mf{k}\subseteq \mf{g}^{(i)}$ for every $i$. Conversely, if $\lie{k}\in \cap_r \: \wmf{g}_r$, then $\lie{k}y=0$. This means that there is an $s$ such that 
$\mf{k}= \mf{g}^{(i)}$ for all $i>s$. $\Box$
\begin{lemma}
    \label{lemma:bifillbasis}
    Let $(M^{i,j})_{i\geq a, j\geq b}$ be a family of vector spaces such that (i) $M^{i,j}\supseteq M^{i,j+1}$, (ii) $M^{i,j}\supseteq M^{i+1,j}$, and (iii) $M_{i+1, j+1}=M^{i,j+1}\cap M^{i+1, j}$. Let $M=M^{a,b}$. Then there is a basis $B=\dot\cup_{i,j} B_{i,j}$ such that:
    \begin{enumerate} 
    \item  $\dot\cup_{i\geq i_0, j\geq j_0} B_{i,j}$ is a basis for $M^{i_0 , j_0}$, and 
    \item $B_{i,j}$ is a basis for $M^{i,j}/(M^{i+1,j}+M^{i,j+1})$.     
    \end{enumerate}
\end{lemma}

\noindent
The proof is straightforward.

We call $B$ above as the filtered basis for the family $(M^{i,j})$. Note that $|B_{i,j]}|=dim(M^{i,j}/(M^{i+1,j}+M^{i,j+1}))$. Also note that if $B'=(B'_{i.j})$ is another filtered basis, then for all $i,j$ and all $b'\in B'_{i,j}$, $b'=b+k+h$ where $b$ is an invertible linear combination of elements of $B_{i,j}$, $k\in M^{i,j+1}$ and $h\in M^{i+1,j}$.

\begin{defn} \label{defn:filbasis}
    The filtered basis $B_{\mf{g}}$ or simply $B$ for $\gamma$ is the basis obtained by applying Lemma \ref{lemma:bifillbasis} to $(\mf{g}^{i,j})$. 
\end{defn}  
 
\begin{prop} \label{lemma:gammahat}
Let $\mf{g}_r$ be the set $\{ \lie{g}|\delta(\lie{g})\geq r \}$. Then we have:
\begin{enumerate}
    \item  $\mf{g}_0 = \mf{g}$ and $\mf{g}_{r}
    \supseteq \mf{g}_{r+1}$. Moreover, there is a smallest $s$ such that $\mf{g}_{s+1}=\mf{k}$. 
    \item Let $\widehat{\mf{g}}_r$ be the subspace generated by the elements $\widehat{\lie{g}}$ with $\lie{g}\in \mf{g}_r$. 
    Then $\widehat{\mf{g}}_r$ is a Lie subalgebra for every $r\geq 0$. Moreover, $\widehat{\mf{g}}_0 = \mf{g}$. 
    \item For every $r$ and $\widehat{\mf{g}}_r \supseteq \widehat{\mf{g}}_{r+1}$ and $\widehat{\mf{g}}_{s+1}=\widehat{\mf{k}}$. 
    \item Finally, $\widehat{\mf{g}}_1 \subseteq  \mf{h}$ and $\widehat{\mf{g}}_2\subseteq \mf{h}_{\overline{y_1}}$. 
\end{enumerate}
\end{prop}

\noindent
{\bf Proof}: (1) follows from the fact that $\delta (\lie{g})\geq 0$ for any $\lie{g}\in \mf{g}$, and Lemma \ref{lemma:bifil}.

Let us prove (2). Let $\widehat{\lie{g}}, \widehat{\lie{g}'}\in \widehat{\mf{g}}_r$ with $\lie{g}, \lie{g}' \in \mf{g}_r$. In order to compute $[\widehat{\lie{g}}, \widehat{\lie{g}'}]$, we examine:
\[ 
[\lie{g}(t), \lie{g}'(t)]y(t)=
\lie{g}(t)\lie{g}' (t)y (t)-\lie{g}'(t) \lie{g}(t) y(t) \\ 
\]
Now:
\[ \begin{array}{rcl} 
\lie{g}(t) (\lie{g}'(t)y(t)) &=& \lie{g}(t) ((\lie{g'}y)(t)) = \lie{g}(t)(\gamma (t) (\lie{g'}y))\\
&=& \lie{g} (t) [t^{d(\lie{g}')+r} y'_{d(\lie{g}')+r} \: + \ldots ]\\
&=& (t^{d(\lie{g})} \lie{g}_{d(\lie{g})} \: + \ldots )[t^{d(\lie{g}')+r} y'_{d(\lie{g}')+r} \: + \ldots ] \\
&=&[t^{d(\lie{g})+d(\lie{g}')+r} y''_{d(\lie{g})+d(\lie{g}')+r} \: + \ldots  ]
\end{array}
\]
Similarly, we have:
\[ \lie{g'}(t) (\lie{g}(t)y(t)) = t^{d(\lie{g})+d(\lie{g}')+r} y'''_{d(\lie{g})+d(\lie{g}')+r} \: + \ldots \]
Thus, we have: 
\[ 
[\lie{g}, \lie{g}'](t) y(t)=[\lie{g}(t), \lie{g}'(t)]y(t)=
t^{d(\lie{g})+d(\lie{g}')+r} v_{d(\lie{g})+d(\lie{g}')+r} \: + \ldots 
\]
for a suitable set $v_i\in V$.
We also have:
\[ [ \lie{g}, \lie{g'}](t)= [ \wl{g}, \wl{g'}]t^{d(\lie{g})+d(\lie{g'})}+\ldots \]
Now, either (i) $d([\lie{g}, \lie{g}'])=d(\lie{g})+d(\lie{g}')$ and $\widehat{[\lie{g}, \lie{g'}]}=[\widehat{\lie{g}}, \widehat{\lie{g}'}]$ 
or (ii) $[\widehat{\lie{g}}, \widehat{\lie{g}'}]=0$. If (i) holds, then the above two equations tells us that $[\lie{g}, \lie{g}']\in \mf{g}_r$. If (ii) holds, then $[\widehat{\lie{g}}, \widehat{\lie{g}'}]=0$. Thus, in both cases $[\widehat{\lie{g}}, \widehat{\lie{g}'}]\in \widehat{\mf{g}}_r$. 

(3) follows from (1) and (2). For (4), note that for any $\lie{g}\in \mf{g}$, if $\lie{g}(t)= g_k t^k+g_{k+1} t^{k+1}+\ldots$, with $k=d(\lie{g})$ and $g_i \in \mf{g}$, and $y(t)=y_0 +y_1 t^1 +\ldots $, then 
\[(\lie{g}y)(t)= (g_k t^k+g_{k+1} t^{k+1}+\ldots)(y_0 +y_1 t^1 +\ldots)\]
Thus if $\lie{g}\in \mf{g}_1$, then $g_k y_0= 0$ and thus $g_k=\widehat{\lie{g}}\in \mf{h}$. If $\lie{g}\in \mf{g}_2$, then we also have $g_k y_1 +g_{k+1} y_0=$, and thus $\wl{g} \in \mf{h}_{\overline{y_1}}$. This proves the assertion. $\Box$

\eat{
Something like the following looks true..
\begin{lemma} We have 
$\mf{g}_0 = \mf{g}$. Moreover $\mf{g}_{r+1}$ may be related to $\widehat{\mf{k}_r(t)}$  where
    $$ \mf{k}_r(t) = \{ \lie{k}(t) \in {\mf{g}}(t) | \lie{k}(t) . (y_0 + y_1t^1 + \ldots + y_r t^{r}) = 0 \}$$
Note that, $\mf{k}_r(t)$ is the Lie algebra stabilizer of the approximation of $y(t)$ upto degree $r$.
Another related space seems to be
$$\{ \lie{g}(t) \in \mf{g}(t) | d(\lie{g}(t) y(t) \geq d(\lie{g}(t)) + r\}$$
Further, $\widehat{\mf{g}}_0 = \mf{g}$ and
$\widehat{\mf{g}}_{r+1} = \{ \lie{g} \in \widehat{\mf{g}}_{r} | \lie{g}.y_r = 0 \} = \{\lie{g} \in \mf{g}| \lie{g}.y_0 = \lie{g}.y_1 = \ldots = \lie{g}.y_r = 0\}$. 
\end{lemma}
}

\vspace*{0.3cm}

\noindent
{\it The algebras $\wmf{g}_r$ will be called the {\bf cancellation algebras}  for the data $(\gamma, y)$.}
\begin{prop}
    \label{prop:Wfilter}
    Let $B=(B_{i,j})$ be the filtered basis as in Defn \ref{defn:filbasis}. Then, we have:
    \begin{enumerate}
        \item The set $B_{*,s+1}=\cup_{i=a}^b B_{i,s+1}$ is a basis for $\mf{k}$.
        \item The set $B'=\cup_{i=a}^b \cup_{j=a}^s$ extends the basis $B_{*,s+1}$ of $\mf{k}$ to that of $\mf{g}$. 
    \end{enumerate}
\end{prop}

\noindent
{\bf Proof}: The first assertion follows from lemma \ref{lemma:bifil} that $\mf{g}^{a,s+1}=\mf{k}$. The second assertion follows from the fact that $\mf{g}^{a,b}=\mf{g}$, and the construction of $B$ as in lemma \ref{lemma:bifillbasis}. $\Box$

\begin{lemma} \label{lemma:gr}
    Let $\lie{g}$ be such that $d(\lie{g})=i$ and $d(\lie{g}y)=j=i+r$. Then 
    \begin{enumerate}
        \item There is a linear combination  $\lie{b}$ of elements of $B_{i,j}$ and a $\lie{k}$ such that $\delta (\lie{k})\geq r+1$ so that $\wl{g}=\wl{b}+\wl{k}$, where either one or both of $\lie{b}$ and $\lie{k}$ being non-zero, and $d(\lie{k})=i$, if $\lie{k}\neq 0$.
        \item There is a linear combination  $\lie{b}$ of elements of $B_{i,j}$ and a $\lie{h}$, so that $\wl{gy}=\wl{by}+\wl{hy}$ and $d(\lie{h}y)=j$, where either one or both of $\lie{b}$ and $\lie{h}$ being non-zero, $\delta (\lie{h})\leq r-1$, if $\lie{h}\neq 0$.
    \end{enumerate}    
\end{lemma}

\noindent
{\bf Proof}: By the hypothesis, we have $\wl{g}\neq 0$ and $\wl{gy}\neq 0$. By the construction of $B$ in \ref{lemma:bifil}, we know that $\C B_{i,j}+ (\mf{g}^{i+1,j}+\mf{g}^{i,j+1})=\mf{g}^{i,j}$. Thus, for (1), we may write:
\[ \lie{g}= \lie{b}+\lie{h}+\lie{k}\]
where $\lie{b}$ is a linear combination of elements in $B_{i,j}$, $\lie{h}\in \mf{g}^{i+1,j}$ and $\lie{k}\in \mf{g}^{i, j+1}$ and $d(\lie{k})=i$. 
Note that if $\lie{g}\in \mf{g}^{i+1,j}+\mf{g}^{i,j+1}$ then $\lie{b}=0$. The condition on $\lie{h}$ implies that $\wl{h}$ does not contribute to $\wl{g}$, while the condition on $\lie{k}$ implies that $\delta (\lie{k})>r$. Hence $\wl{g}=\wl{b}+\wl{k}$, with $\lie{k}$ having the stated properties. This proves (1). 

For (2), 
we may write:
\[ \lie{g}= \lie{b}+\lie{h}+\lie{k}\]
where $\lie{b}$ is a linear combination of elements in $B_{i,j}$, $\lie{h}\in \mf{g}^{i+1,j}$ and $\lie{k}\in \mf{g}^{i, j+1}$ and $d(\lie{h}y)=j$. These two conditions imply that $\wl{ky}$ does not contribute to $\wl{gy}$, hence $\wl{gy}=\wl{by}+\wl{hy}$. Furthermore, since $d(\lie{h})\geq i+1$, we have $\delta (\lie{h}) = j - d(\lie{h}) \leq  r - 1$. This proves (2) and the lemma. $\Box$ 

\begin{prop} \label{prop:hatbasis}
    Let $B$ be as above. Then:
    \begin{enumerate}
        \item For every $i,j$, let $C_{i,j}= \{ \widehat{\lie{b}}| \lie{b}\in B_{i', j'} \: with \: i'\geq i, j'\geq j\}$. Then $C_{i.j}$ is a basis for the vector space $\widehat{\mf{g}^{i,j}}$. 
         \item For every $i,j\leq s$, let $D_{i,j}= \{ \widehat{\lie{b}y}| \lie{b}\in B_{i', j'} \: with \: i'\geq i, s\geq j'\geq j\}$. Then $D_{i.j}$ is a basis for the vector space $\widehat{\mf{g}^{i,j}y}$. 
    \end{enumerate}
\end{prop}

\noindent
{\bf Proof}: We begin with (1) and prove this by reverse induction. Let us assume the assertion to be true for all $i',j'$ where either $i'>i$ or $j'>j$, and let us prove it for $i,j$ and let $r=j-i$. Let us pick a $\lie{g}\in \mf{g}^{i,j}$ such that $\wl{g}\neq 0$. Then by lemma \ref{lemma:gr} (1), we have a  $\lie{b}$, linear combination of elements of $B_{i,j}$ and a $\lie{k}$ such that $\delta (\lie{k})\geq r+1$ so that $\wl{g}=\wl{b}+\wl{k}$ and $d(\lie{k})=i$, if $\lie{k}\neq 0$. By induction, $\wl{k}\in \C C_{i,j+1}$ and thus the claim holds for $C_{i,j}$ and (1) is proved. The proof of (2) is similar. $\Box$

\begin{corollary} \label{cor:grseries}
 We have:
\begin{enumerate}
    \item For any element $\lie{g}\in \mf{g}_r$ such that $d(\lie{g})=i$, there is a sequence $i+r\leq j_1 < \ldots < j_s$ such that $\wl{g}=\wl{b}_1+\ldots +\wl{b}_s$ with $\lie{b}_k \in \C B_{i,j_k}$. 
\item For any element $\lie{g}\in \mf{g}_r$, such that $j=d(\lie{g}y)$, there is a sequence $j-r\leq i_1 < \ldots < i_s$ such that $\wl{gy}=\wl{by}_1+\ldots +\wl{by}_s$, with each $\wl{b}_k \in \C B_{i_k,j}$.
\end{enumerate}
\end{corollary}

\subsection{The cancellation algebras and filtered tangent spaces.}

\begin{prop}
    \label{prop:gr1}
    For the basis $B$ as above. Let $\widehat{B}_{i,j}=\{ \widehat{\lie{b}}| \lie{b}\in B_{i,j}\}$ and $\widehat{B}=\{ \wl{b}| \lie{b}\in B \}$. Then 
     $\{ \widehat{B}_{i,j}|j-i\geq r \}$ is a basis for $\wmf{g}_r$. Moreover $\wmf{g}_r /\wmf{g}_{r+1}$ has a basis consisting of cosets of $\dot \cup_{j-i=r} \widehat{B}_{i,j}$. 
\end{prop}

\noindent
{\bf Proof}: 
The proof is by downward induction. 
We start at $r=D$, where $\mf{g}_D = \mf{k}$. The basis $B_{*,s+1}$ is included at every level $r \leq D$ since $d(\lie{k}y) = \infty$ for every element $\lie{k} \in \mf{g}_D $.

Let us assume the statement to be true for $r+1$ and let us prove it for $r$. We know that 
$\wmf{g}_r$ is the space generated by $C_r=\{ \wl{g}|\delta (\wl{g})=r\}$ and $ \wmf{g}_{r+1}$. 
Suppose $\lie{g}\in C_r$ is such that $d(\lie{g})=i$ and $d(\lie{g}y)=i+r$. Then by lemma \ref{lemma:gr} (1), there is a linear combination  $\lie{b}$ of elements of $B_{i,j}$ and a $\lie{k}\in \mf{g}_{r+1}$ so that $\wl{g}=\wl{b}+\wl{k}$.
Going modulo $\wmf{g}_{r+1}$, we have $\wl{g} + \wmf{g}_{r+1} = \wl{b}+\wmf{g}_{r+1}$, so this element is in the span of cosets of $\dot \cup_{j-i=r} \widehat{B}_{i,j}$. Independence 
follows from the independence of elements $B_{i,j}$, with $j - i \geq r$.
$\Box$

 \begin{prop} \label{prop:gry}
With the above notation, we have:
\begin{enumerate}
    \item 
Let $\C B'y\subseteq V$ be the vector space generated by $B'y =\{ \lie{b}y| \lie{b}\in B_{i,j} \: with \: j\leq s\}$. Then $\C B'y =T_y$, the tangent space $T_y O(y)$, and $B'y$ is a basis.
\item Let $\widehat{Ty}$ be the vector space generated by the leading terms of $T_y$, then $dim(\widehat{Ty})=dim(T_y)$ and $\widehat{B'y}=\{ \wl{by}| b\in B' \}$ is a basis.
\item  Let $Q_r \subseteq V$ be the space generated by $\{ \widehat{\lie{g}y}| \delta(\lie{g})\leq r \}$. Then 
       $\{ \wl{by}|\lie{b}\in B_{i,j} \: with \: j-i\leq r \}$ is a basis for $Q_r$. Moreover $Q_r /Q_{r-1}$ is generated by  $\dot\cup_{j=i+r} \widehat{By}_{i,j}$, where $\widehat{By}_{i,j}=\{ \wl{by}|\lie{b}\in B_{i,j} \}$. 
 \end{enumerate}
 \end{prop}

\noindent
{\bf Proof}: By Prop. \ref{prop:Wfilter} (2), $B'=\dot\cup_{j\leq s} B_{i,j}$ is a basis for $\mf{g}/\mf{k}$, and we also have $\mf{k}y=0$. Thus $B'y$ is a basis for $T_y$. This proves (1). 
For (2), ampleness implies the first part. For the second part, since $\mf{g}^{a,a}=\mf{g}$, by Prop. \ref{prop:hatbasis}(2), $D_{a,a}$ is a basis for $\widehat{\mf{g}y}$. But $D_{a,a}$ is precisely $\widehat{B'y}$. 

We prove (3) by induction on $r$ starting with $r=0$. Clearly $Q_{-1} = 0$, so the base case holds. Assume $r \geq 1$ and the statement holds for $r-1$. 
Let $\lie{g}$ be such that $d(\lie{g})=i$ and $\delta (\lie{g})=r$, i.e., $d(\lie{g}y)= j = i+r$. 
Then by lemma \ref{lemma:gr} (2), 
there is a linear combination  $\lie{b}$ of elements of $B_{i,j}$ and a $\lie{h}$ with
$\delta(\lie{h}) \leq r-1$, so that $\wl{gy}=\wl{by}+\wl{hy}$. The first term is generated by elements in $\widehat{By}_{i,j}$ and the second term is in $Q_{r-1}$. Independence follows from (1). This proves (3) and the proposition. $\Box$

\begin{corollary} \label{cor:dimequal}
(i) For every $r$, $dim(Q_r / Q_{r-1})=dim(\wmf{g}_r/ \wmf{g}_{r+1})=  \sum_i |B_{i,i+r}|$. (ii) Moreover, $Q_0 =T_z O(z)$ and $\wmf{g}_1=\mf{h}$.   
\end{corollary}

\noindent
{\bf Proof}: The first part is clear since both $Q_r/Q_{r-1}$ and $\wmf{g}_r/ \wmf{g}_{r+1}$ are indexed by the same basis. Coming to (ii), we already know that $\wmf{g}_1\subseteq \mf{h}$. 
Let $\lie{h}\in \mf{h}-\wmf{g}_1$. Since $\widehat{B}$ is a basis for $\mf{g}=\wmf{g}_0$, we may express $\lie{h}=\widehat{b}_{i_1,j_1} +\ldots +\widehat{b}_{i_s , j_s}$, where $b_{i_k, j_k} \in \C B_{i_k , j_k}$. If $j_k-i_k>0$, then $\widehat{b}_{i_k , j_k}\in \wmf{g}_1\subseteq \mf{h}$. Thus, we may assume that $i_k=j_k$ for the elements $\widehat{b}_{i_k, j_k}$. This implies that $\widehat{b_{i_k, j_k}y}= \widehat{b_{i_k , j_k}}\widehat{y}=\widehat{b_{i_k , j_k}}z$.
Now, the condition $\lie{h}z=0$ implies that $ \widehat{b_{i_1,j_1}y} +\ldots +\widehat{b_{i_s , j_s}y}=0$. This contradicts that such elements form a basis of $Q_0$. This proves $\wmf{g}_1=\mf{h}$ and consequently $Q_0 =T_z$. $\Box$

\vspace*{0.3cm}

{\it We conjecture that the flag $(Q_r)$ provides information about the singular locus of $O(z)$ within $\overline{O(y)}$.} 

\vspace*{0.3cm}

\noindent
The following lemma outlines the action of $\widehat{\mf{g}}_s$ on $(\widehat{T_y})_r$ in terms of the action of $\mf{g}$ on $T_y$.

\begin{lemma}
    \label{lemma:rcomput}
    Let $\lie{h}\in \mf{g}_s$ such that $\lie{h}(t)=h_a t^a +\ldots$, i.e., $d(\lie{h})=a$ and $d(\lie{h}y)=a+s$. Let $r<s$ and $\lie{g}\in \mf{g}_r$ such that $d(\lie{g})=b$ while $deg(\lie{g}y)=b+r$. 
    Let $\lie{g}'=[\lie{h}, \lie{g}]$. Then (a) $\sum_{i=0}^{c} g'_{a+b+i} y_{c-i}=0$ for all $c<r$ and (b) $\widehat{\lie{h}}\cdot \widehat{\lie{g}y}=\sum_{i=0}^r g'_{a+b+i} y_{r-i}$. 
We denote the last term by $(\widehat{\lie{g}' y})_{a+b+r}$.
\end{lemma}

\eat{\noindent
{\bf Proof}:  This is a long computation, which prove in detail for $a=b=0$, the general case being a simple extension.  Note that $\lie{h}\in \mf{g}_s$ implies that (i) for all $i< s$ we have $\sum_{j=0}^i h_j y_{i-j}=0$. Since $\widehat{\lie{g}y}=t^r y'_r +\ldots$, we have (ii) for all $j<r$ we have $\sum_{i=0}^j g_i y_{j-i} =0$. 

Now coming to the assertion, let $c\leq r$. We have:
\[ \begin{array}{rcl}
h_0 (g_0 y_c +\ldots + g_c y_0 )&=& \sum_{i=0}^c [h_0 , g_i]y_{c-i} +\sum_{i=0}^c g_i h_0 y_{c-i} \\
&=& \sum_{i=0}^c [h_0 , g_i]y_{c-i} -\sum_{i=0}^c g_i (\sum_{j=1}^{c-i} h_j y_{c-i-j}) \mbox{    (by applying  (i) and with $c<s$)}\\
&=& \sum_{i=0}^c [h_0 , g_i] y_{c-i}+\sum_{i=0}^c \sum_{j=1}^{c-i} [h_j, g_i]y_{c-i-j} -  \sum_{i=0}^c \sum_{j=1}^{c-i} h_j g_i y_{c-i-j}\\
&=& \sum_{i=0}^c [h_0 , g_i] y_{c-i}+\sum_{i=0}^c \sum_{j=1}^{c-i} [h_j, g_i]y_{c-i-j} -  \sum_{j=1}^c h_j (\sum_{i=0}^{c-j} g_i y_{c-i-j})\\
&=& \sum_{i=0}^c [h_0 , g_i] y_{c-i}+\sum_{i=0}^c \sum_{j=1}^{c-i} [h_j, g_i]y_{c-i-j} \mbox{   (by applying (ii))}\\
&=& \sum_{i=0}^c (\sum_{j=0}^i ([h_j , g_{i-j}])y_{c-i}\\
&=& \sum_{i=0}^c g'_i y_{c-i}
\end{array} \]
Now for $c<r$, we have $(g_0 y_c +\ldots + g_c y_0 )=0$ and hence:
\[0=h_0 (g_0 y_r +\ldots + g_r y_0 )=\sum_{i=0}^c g'_i y_{c-i} \]
This proves the first assertion (a), with $a=b=0$. The general case is achieved by managing indices.  

When $c=r$, we have:
\[\widehat{\lie{h}}\cdot \widehat{\lie{g}y}=h_0 (g_0 y_r +\ldots + g_r y_0 )=\sum_{i=0}^r g'_i y_{r-i} \]
This proves (b). $\Box$

\textcolor{red}{This is what I mean - I still don;t see anything wrong}
{\bf Proof}
}
\noindent
{\bf Proof}:
We compute $\gamma(t)(\lie{g'}y)$ in two ways. 
\eat{Assume
\[
d(h)=a,\qquad d(g)=b,\qquad
d(hy)=a+s,\qquad d(gy)=b+r,
\qquad s\ge r+1.
\]
}
Since $\lie{g'}(t)=[\lie{h}(t),\lie{g}(t)]$, we have $d(\lie{g'})\ge a+b$.
Let us write
\[
\lie{g'}(t)=\sum_{i\ge0}t^{a+b+i}g'_{a+b+i},
\qquad
y(t)=\sum_{j\geq 0}t^jy_j.
\]
The coefficient of $t^{a+b+c}$ in
$\gamma(t)(\lie{g'}y)=\lie{g'}(t)y(t)$ is therefore
\[
\sum_{i=0}^c g'_{a+b+i}y_{c-i}.
\]

On the other hand,
\[
\gamma(t)(\lie{g'}y)
=\lie{h}(t)\bigl(\lie{g}(t)y(t)\bigr)
-\lie{g}(t)\bigl(\lie{h}(t)y(t)\bigr).
\]
Using
\[
\begin{aligned}
\lie{h}(t)&=t^a\widehat{\lie{h}}+t^{a+1} h_{a+1} + \ldots,\\
\lie{g}(t)y(t)&=t^{b+r}\widehat{\lie{g}y}+t^{b+r+1} y'_{b+r+1} + \ldots,
\end{aligned}
\]
the first term equals
\[
t^{a+b+r}\widehat{\lie{h}}\cdot\widehat{\lie{g}y}
+O(t^{a+b+r+1}).
\]
Since $d(\lie{g}) = b$, and $d(\lie{h}y) = a+s$, the second term has order at least $a+b+s\ge a+b+r+1$.
Consequently,
\[
\gamma(t)(\lie{g'}y)
=t^{a+b+r}\widehat{\lie{h}}\cdot\widehat{\lie{g}y}
+O(t^{a+b+r+1}).
\]

Comparing coefficients we get 
\[
\begin{array}{rcl}

\sum_{i=0}^c g'_{a+b+i}y_{c-i}&=&0
\qquad(0\le c<r)\\

\sum_{i=0}^r g'_{a+b+i}y_{r-i}
&=&\widehat{\lie{h}}\cdot\widehat{\lie{g}y}.
\end{array} \]
This proves the proposition. $\Box $

\begin{corollary} \label{cor:liebr}
    Let $\lie{h}, \lie{g}$ be such that $\delta(\lie{h})=\delta(\lie{g})=r$ and let $\lie{g'}=[\lie{h}, \lie{g}]$. Let $d(\lie{h})=a$ and $d(\lie{g})=b$ and $\lie{g'}(t)=\sum_{i\geq 0}g'_{a+b+i}t^{a+b+i}$. Then (i) for $c<r$, we have $\sum_{i=0}^c g'_{a+b+i}y_{c-i}=0$, and (ii) $
\wl{h}\cdot\widehat{\lie{g}y}
-\wl{g}\cdot\widehat{\lie{h}y}
=
\sum_{i=0}^{r}g'_{a+b+i}y_{r-i}
=
[t^{a+b+r}]\,\gamma(t)(\lie{g}'y)$.  In particular, $d(\lie{g}'y)\geq a+b+r$.
If equality holds, the expression in (ii) equals
$\widehat{\lie{g}'y}$; if the inequality is strict,
that expression is zero.
\end{corollary}

\noindent
{\bf Proof}: The first part follows from the earlier computation. For (ii)
we calculate,
\[
\gamma(t)(\lie{g}'y)
=
\lie{h}(t)\bigl(\gamma(t)(\lie{g}y)\bigr)
-
\lie{g}(t)\bigl(\gamma(t)(\lie{h}y)\bigr).
\]
The hypotheses give
\[
\begin{aligned}
\lie{h}(t)&=t^a\wl{h}+O(t^{a+1}),&
\gamma(t)(\lie{h}y)
 &=t^{a+r}\widehat{\lie{h}y}+O(t^{a+r+1}),\\
\lie{g}(t)&=t^b\wl{g}+O(t^{b+1}),&
\gamma(t)(\lie{g}y)
 &=t^{b+r}\widehat{\lie{g}y}+O(t^{b+r+1}).
\end{aligned}
\]
Consequently,
\[
\gamma(t)(\lie{g}'y)
=
t^{a+b+r}
\left(
\wl{h}\cdot\widehat{\lie{g}y}
-\wl{g}\cdot\widehat{\lie{h}y}
\right)
+O(t^{a+b+r+1}).
\]
Comparing this with
\[
\lie{g}'(t)y(t)
=
\sum_{c\geq0}t^{a+b+c}
\left(\sum_{i=0}^{c}g'_{a+b+i}y_{c-i}\right)
\]
proves the assertions.
$\Box$

\begin{lemma}
    \label{lemma:4cases}
    Let $\lie{h}\in \mf{g}_s$ such that $d(\lie{h})=a$ and $d(\lie{h}y)=a+s$. Let $r<s$ and $\lie{g}\in \mf{g}_r$ such that $d(\lie{g})=b$ while $deg(\lie{g}y)=b+r$. 
    Let $\lie{g}'=[\lie{h}, \lie{g}]$. Then four cases emerge:
    \begin{enumerate}
        \item[(i)] 
    $d(\lie{g}')= a+b$ and $d(\lie{g}'y)=a+b+r$, and then $\widehat{\lie{h}}\cdot \widehat{\lie{g}y}=\widehat{\lie{g}'y}$, $[\wl{h}, \wl{g}]=\wl{g'}$ and $\delta(\lie{g}') = r$.
    \item[(ii)] $d(\lie{g}')= a+b$ and $d(\lie{g}'y)>a+b+r$ and then $\widehat{\lie{h}}\cdot \widehat{\lie{g}y}=0, [\wl{h}, \wl{g}]=\wl{g'}$ and $\delta(\lie{g}')\geq r+1$ 
    \item[(iii)] $d(\lie{g}')> a+b$ and $d(\lie{g}'y)>a+b+r$ and then $ \widehat{\lie{h}}\cdot \widehat{\lie{g}y}=0$ and $[\wl{h}, \wl{g}]=0$.
    \item[(iv)] $d(\lie{g}')> a+b$ and $d(\lie{g}'y)=a+b+r$ and then $ \widehat{\lie{h}}\cdot \widehat{\lie{g}y}=\widehat{\lie{g'}y} \neq 0$, $\delta(\lie{g'})\leq r-1 $ and $[\wl{h}, \wl{g}]=0$.
     \end{enumerate}
\end{lemma}

\noindent
{\bf Proof}: Lemma \ref{lemma:rcomput} shows that  $d(\lie{g}'y) \geq a+b +r$ and $\widehat{\lie{h}}\cdot \widehat{\lie{g}y}=(\widehat{\lie{g}' y})_{a+b+r}$. The four cases are then obvious. $\Box$

\begin{prop}
\label{prop:Tbijection}
    Let $Q_r \subseteq V$ be the space generated by $\{ \widehat{\lie{g}y}| \delta(\lie{g})\leq r \}$, as in Prop. \ref{prop:gry}. Then:
    \begin{enumerate}
        \item[(i)] $Q_r$ is a $\widehat{\mf{g}}_{r+1}$-module.  
        \item[(ii)]  $Q_r / Q_{r-1} \cong \widehat{\mf{g}}_r/\widehat{\mf{g}}_{r+1}$ as $\widehat{\mf{g}}_{r+1}$-modules. 
    \end{enumerate}
\end{prop}

\noindent
{\bf Proof}: For (i), we know that by Prop. \ref{prop:Wfilter} the sets $\widehat{By}_{\leq r}=\{ \widehat{By}_{i,j} |j-i\leq r \}$ and  $\widehat{By}_{>r}=\{ \widehat{By}_{i,j} |j-i > r \}$
are linearly independent, and that $\widehat{B}_{>r}$ is a basis for $\widehat{\mf{g}}_{r+1}$.
 The space $Q_r\subseteq V$ and $V$ itself is an $\mf{g}$-module. Thus to prove (i), we need to merely show that for every $\lie{b}\in \widehat{By}_{\leq r}$ and $\lie{b'}\in \widehat{B}_{> r}$, we have $ \widehat{\lie{b}'}\cdot \widehat{\lie{b}y} \in Q_r$. Now, by Lemma \ref{lemma:4cases}, we see that $ \widehat{\lie{b}'}\cdot \widehat{\lie{b}y}\neq 0$ may happen only in cases (i) and (iv). Let $\lie{g'}=[\lie{b'}, \lie{b}]$.
 For (i) we have $\widehat{\lie{b}'}\cdot \widehat{\lie{b}y}=\widehat{\lie{g'}y}$, $d(\lie{g'})=d(\lie{b}')+d(\lie{b})$ and $\delta (\lie{g'})=r$.
 This implies that  $\widehat{\lie{g'}y} \in Q_r$. For case (iv), we see that $\delta (\lie{g'})\leq r-1$, and thus $\widehat{\lie{g'}y} \in Q_{r-1}$. This proves (i).

Coming to (ii), let $S_r$ be the module $\widehat{\mf{g}}_r / \widehat{\mf{g}}_{r+1}$. The set $\{ \widehat{\lie{b} y}+Q_{r-1}| \lie{b} \in B_{ r}\}$ is a basis for $Q_r / Q_{r-1}$, while $\{ b +\widehat{\mf{g}}_{r+1}| b \in B_r\}$ is a basis for $S_r$.  
Thus there is a vector space isomorphism $\phi:  S_r \rightarrow Q_r $ given by $\phi (\lie{b} +\widehat{\mf{g}}_{r+1})=\widehat{\lie{b}y}+Q_{r-1}$. 
For any $\lie{b}\in B_{\leq r}$ and $\lie{b'}\in B_{> r}$ and let $\lie{g}'=[\lie{b}', \lie{b}]$. Let $d(\lie{b}')=a$ and $d(\lie{b})=b$. Let $d(\lie{b}y)=b+r$. We know that $\lie{g}' \in \mf{g}^{a+b}\cap \mf{g}^{(a+b+r)}$. 

By lemma \ref{lemma:4cases}, we have the four cases:
\begin{enumerate}
    \item 
    $d(\lie{g}')= a+b$,  $d(\lie{g}'y)=a+b+r$, and $\widehat{\lie{b}'}\cdot \widehat{\lie{b}y}=\widehat{\lie{g}'y}$. In this case $\delta (\lie{g}'
    )=r$ and thus $\widehat{\lie{g}'}=[\widehat{\lie{b}'}, \widehat{\lie{b}}] \in \wmf{g}_r$. 
    By Cor. \ref{cor:grseries} (1), we then have $\widehat{\lie{g}'}= [\widehat{\lie{b}'}, \widehat{\lie{b}}]+ \sum_{j>a+b+r} \wl{b}_{j}$, with $\lie{b}_j\in \C B_{a+b, j}$ and $[\wl{b'}, \wl{b}]\in \C \widehat{B}_{a+b, a+b+r}$. 
    This also implies that $\widehat{\lie{g}'y}= \widehat{[\widehat{\lie{b}'}, \widehat{\lie{b}}]y} $.
    Thus the action of $\lie{b}'$ on $Q_r$ matches that of $S_r$. 
    \item $d(\lie{g}')= a+b$ and $d(\lie{g}'y)>a+b+r$ and then $\widehat{\lie{b}'}\cdot \widehat{\lie{b}y}=0$ and $\lie{g}'\in \mf{g}_{r+1}$. Thus the action of $\lie{b}'$ is zero on both $Q_r$ and $S_r$.  
    \item $d(\lie{g}')> a+b$ and $d(\lie{g'}y)>a+b+r$, then $ \widehat{\lie{b}'}\cdot \widehat{\lie{b}y}=0$ and $[\widehat{\lie{b'}}, \widehat{\lie{b}}]=0$. Once again, the action of $\lie{b}'$ is zero on both $Q_r$ and $S_r$.  
    \item $d(\lie{g}')> a+b$ and $d(\lie{g'}y)=a+b+r$, and then $ \widehat{\lie{b}'}\cdot \widehat{\lie{b}y}=\widehat{\lie{g'}y}$ with $\delta (\lie{g'})\leq r-1$. This implies that $\widehat{\lie{b}'}\cdot \widehat{\lie{b}y}=0+Q_{r-1}$. This matches with $[\widehat{\lie{b'}}, \widehat{\lie{b}}]=0$. Once again, the action of $\lie{b}'$ is zero on both $Q_r$ and $S_r$. 
\end{enumerate}
The above proves (iv) and the proposition. $\Box$

\begin{prop} \label{prop:derivation}
     Let $D_r : \wmf{g}_r \rightarrow V/Q_{r-1}$ be defined as follows. Let $d(\lie{g})=a$ with $\lie{g}(t)= t^a g_a +\ldots $. We define $D_r (\lie{g})= \sum_{i=0}^r g_{a+i}y_{r-i}+ Q_{r-1}$. Then $D_r:\wmf{g}_r \rightarrow Q_r/Q_{r-1}$ is well defined and is a derivation. Moreover $D(\wmf{g}_{r+1})=0$. 
    
\end{prop}
\noindent
{\bf Proof}:
We see that $D_r$ is merely the composition
\[
\wmf{g_r}\longrightarrow\wmf{g_r}/\wmf{g_{r+1}}
\xrightarrow{\phi_r}Q_r/Q_{r-1}
\lhook\joinrel\longrightarrow V/Q_{r-1}.
\]
Moreover, that it is a derivation follows from Corollary \ref{cor:liebr}. $\Box$

\eat{\begin{lemma}
    \label{lemma:deriv}
    Let $\lie{h}\in \mf{g}_r$ such that $d(\lie{h})=a$ and $d(\lie{h}y)=a+r$. Let $\lie{g}\in \mf{g}_r$ such that $d(\lie{g})=b$ while $deg(\lie{g}y)=b+r$. 
    Let $\lie{g}'=[\lie{h}, \lie{g}]$. Then (a) $\sum_{i=0}^{c} g'_{a+b+i} y_{c-i}=0$ for all $c<r$ and (b) $\widehat{\lie{h}}\cdot \widehat{\lie{g}y}=\sum_{i=0}^r g'_{a+b+i} y_{r-i}$. 
We denote the last term by $(\widehat{\lie{g}' y})_{a+b+r}$.
\end{lemma}

\noindent
{\bf Proof}:  As before, we prove this for $a=b=0$, the general case being a simple extension.  Note that $\lie{h}\in \mf{g}_r$ implies that (i) for all $i< r$ we have $\sum_{j=0}^i h_j y_{i-j}=0$. Since $\widehat{\lie{g}y}=t^r y'_r +\ldots$, we have (ii) for all $j<r$ we have $\sum_{i=0}^j g_i y_{j-i} =0$. 

Now coming to the assertion, let $c <r$. We have:
\[ \begin{array}{rcl}
h_0 (g_0 y_c +\ldots + g_c y_0 )&=& \sum_{i=0}^c [h_0 , g_i]y_{c-i} +\sum_{i=0}^c g_i h_0 y_{c-i} \\
&=& \sum_{i=0}^c [h_0 , g_i]y_{c-i} -\sum_{i=0}^c g_i (\sum_{j=1}^{c-i} h_j y_{c-i-j}) \mbox{    (by applying  (i) and with $c<r$)}\\
&=& \sum_{i=0}^c [h_0 , g_i] y_{c-i}+\sum_{i=0}^c \sum_{j=1}^{c-i} [h_j, g_i]y_{c-i-j} -  \sum_{i=0}^c \sum_{j=1}^{c-i} h_j g_i y_{c-i-j}\\
&=& \sum_{i=0}^c [h_0 , g_i] y_{c-i}+\sum_{i=0}^c \sum_{j=1}^{c-i} [h_j, g_i]y_{c-i-j} -  \sum_{j=1}^c h_j (\sum_{i=0}^{c-j} g_i y_{c-i-j})\\
&=& \sum_{i=0}^c [h_0 , g_i] y_{c-i}+\sum_{i=0}^c \sum_{j=1}^{c-i} [h_j, g_i]y_{c-i-j} \mbox{   (by applying (ii))}\\
&=& \sum_{i=0}^c (\sum_{j=0}^i ([h_j , g_{i-j}])y_{c-i}\\
&=& \sum_{i=0}^c g'_i y_{c-i}
\end{array} \]
Now for $c<r$, we have $(g_0 y_c +\ldots + g_c y_0 )=0$ and hence:
\[0=h_0 (g_0 y_r +\ldots + g_r y_0 )=\sum_{i=0}^c g'_i y_{c-i} \]
This proves the first assertion (a), with $a=b=0$. The general case is achieved by managing indices.  

When $c=r$, we have:
\[ \begin{array}{rcl}
h_0 (g_0 y_r +\ldots + g_r y_0 )&=& \sum_{i=1}^r [h_0 , g_i]y_{r-i} +\sum_{i=1}^r g_i h_0 y_{r-i} + h_0 g_0 y_r \\
&=& \sum_{i=1}^r [h_0 , g_i]y_{r-i} -\sum_{i=1}^r g_i (\sum_{j=1}^{r-i} h_j y_{r-i-j}) + h_0 g_0 y_r \mbox{    (by applying  (i))}\\
&=& \sum_{i=1}^r [h_0 , g_i] y_{r-i}+\sum_{i=1}^r \sum_{j=1}^{r-i} [h_j, g_i]y_{r-i-j} -  \sum_{i=1}^r \sum_{j=1}^{r-i} h_j g_i y_{r-i-j} + h_0 g_0 y_r\\
&=& \sum_{i=1}^r [h_0 , g_i] y_{r-i}+\sum_{i=1}^r \sum_{j=1}^{r-i} [h_j, g_i]y_{r-i-j} -  \sum_{j=1}^{r-1} h_j (\sum_{i=1}^{r-j} g_i y_{r-i-j}) +h_0 g_0 y_r\\
&=& \sum_{i=1}^r [h_0 , g_i] y_{r-i}+\sum_{i=1}^r \sum_{j=1}^{r-i} [h_j, g_i]y_{r-i-j} + h_0 g_0 y_r \mbox{(by applying (ii))}\\
&=& \sum_{i=1}^r (\sum_{j=0}^i [h_j , g_{i-j}])y_{r-i}\\
&=& \sum_{i=1}^r g'_i y_{r-i} - \sum_{i=0}^r [h_i, g_0] y_{r-i}  + h_0 g_0 y_r
\end{array} \]
 $\Box$}

\section{The groups $H^r$}
Let $V'=V\otimes {\cal O}$ and $V[r]=V' / t^{r+1}V'$ as defined earlier. Note that $dim(V[r])=(r+1)dim(V)$.
Let $G'=G({\cal O})$, the group of matrices satisfying the equations of $G$ with matrix entries in ${\cal O}$. Note that $V'$ and $t^{r+1}V'$ are $G'$-modules, and so $V[r]$ too is a $G'$-module via the map $\pi_\rho^r : G'\rightarrow GL(V[r])$. 

Let $\C[r]$ be the algebra $\C[t]/(t^{r+1})$ and $G^r $ be the subgroup of matrices in $GL(X)$ with entries in $\C[r]$ satisfying the algebraic relations defining $G$. We then have the natural map $\pi^r :G'\rightarrow G^r$. 

\begin{lemma}
    $G^r$ is an algebraic group and equals the image $\pi^r (G')\subseteq GL(V[r])$. 
\end{lemma}

\begin{prop} \label{prop:jet}
    Let $G^r$ be as above. Then:
    \begin{enumerate}
        \item We have $ker(\pi^r_\rho)\supseteq ker(\pi^r)$, in other words, the map $\pi^r_\rho$ factors through $\pi^r$.
        \item There is an action $\rho^r : G^r \times V[r]\rightarrow V[r]$ derived from $\rho$. This matches the action of $G^r$ through $\pi^r_\rho$. 
        \item The Lie algebra of $G^r$ is 
        \[ \mf{g}\otimes \C[r]= \{ \lie{g}_0 +t\lie{g}_1 +\ldots t^r \lie{g}_r| \lie{g}_i \in \mf{g} \}\]
         This element is denoted as $\lie{g}(t)$.
         \item Moreover, the action of $\lie{g}(t)$ of $v(t)=v_0 +tv_1 +\ldots +t^r v_r \in V[r]$ is given by:
         \[ \lie{g}(t) v(t)= \sum_{i=0}^r t^r \left( \sum_{j=0}^{i} \lie{g}_j v_{i-j} \right)\]
        \item For $s>r$, there are projection operators $\pi^r :G^s \rightarrow G^r, \pi^r : Lie(G^s)\rightarrow Lie(G^r)$ and $\pi^r: V[s]\rightarrow V[r]$ are such that the group and Lie algebra actions are equivariant. 
        
    \end{enumerate}
\end{prop}

\noindent
The proofs are standard.

\eat{{\bf Proof: }
We prove (2) and (4).
For (2) let $f$ be a polynomial equation in the matrix entries that elements of $G$ satisfy. Consider a tangent vector at the identity element of $G^r$,
$I + \epsilon(X_0 + t X_1 + \ldots + t^r X_r)$. Evaluating $f$ at the tangent vector in the ring $\C[t, \epsilon]/(t^{r+1}, \epsilon^2)$ we get,
\[f(I + \epsilon (X_0 + t X_1 + \ldots + t^r X_r)) = f(I) + \epsilon(df_{I}(X_0) + t df_I(X_1) + \ldots + t^r df_I(X_r)).\]
Now $f(I)=0$. The above equation vanishes precisely when $df_I(X_i) = 0$, for all $i$. Since this should hold for all $f$, we conclude that $X_i \in \mf{g}$ for all $i$, proving (2).

\noindent Since $G$ is smooth every jet of $G^r$ lifts to a jet of $G^{r+1}$ proving (4),see ~\cite{stacks-047N,stacks-02H6}. 
$\Box$
}

Let us treat $y(t)=\gamma(t)y$ as an element of $V'$ and let $y^r =\pi^r (y(t))\in V[r]$. 
\begin{prop}    
 Let $H^{r+1}=G^r_{y^r}$ be the stabilizer of $y^r$. Then for the sequence $ y^{s+1} \stackrel{\pi^{s-1}}{\rightarrow} y^{s} {\rightarrow} \ldots y^{i} \stackrel{\pi^{i-1}}{\rightarrow} y^{i-1} \rightarrow \ldots \rightarrow y^1 \stackrel{\pi^{0}}{\rightarrow} y^0$
 we have the chain $ H^{s+1} \stackrel{\pi^{s-1}}{\rightarrow} H^{s} {\rightarrow} \ldots H^{i+1} \stackrel{\pi^{i-1}}{\rightarrow} H^i \rightarrow \ldots \rightarrow H^2 \stackrel{\pi^{0}}{\rightarrow} H^1$, with $H^1=H=G_z$. Moreover, there is a corresponding sequence of algebraic groups $H^{s+1} (0)\subseteq \ldots H^2 (0)\subseteq H^1 (0)=H$. 
\end{prop}

\noindent
{\bf Proof}: That $\pi^{r}(H^{r+2})\subseteq H^{r+1}$ follows from equivariance of $\pi^r$. For the second assertion, note that $V[0]=V$ and $y^0=z\in V$. The final assertion is obtained by successive application of $\pi^i$ and corresponds to evaluation at $t=0$. $\Box$

\begin{prop}
    \label{prop:lie}
    For an element $\lie{b}\in B$, let $\tilde{\lie{b}}(t)=t^{-d(\lie{b})}\lie{b}(t)$. 
    \begin{enumerate}
        \item Then 
\[ \mf{h}^{r+1} (t)=Lie(H^{r+1}) =\{ \lie{g}(t)=\lie{g}_0 +t^1 \lie{g}_1 + \ldots +t^r \lie{g}_r | \lie{g}(t)y^r (t)=0 \} \]
 \item
For $\lie{b}\in B$, put $\nu_{\lie{b}}=\delta(\lie{b})$,
with $\nu_{\lie{b}}=\infty$ when $\lie{b}\in\mf{k}$.
A $\mathbb C$-basis of $\mf{h}^{r+1}(t)=\Lie(H^{r+1})$ is
\[
\left\{
t^a\tilde{\lie{b}}(t)\bmod t^{r+1}
\;\middle|\;
\lie{b}\in B,\quad
\max\{0,r+1-\nu_{\lie{b}}\}\leq a\leq r
\right\}.
\]
For $\nu_{\lie{b}}=\infty$, the lower bound is understood
to be $0$. An empty range contributes no basis elements.
 \eat{ \item A $\C$-basis for $\mf{h}^{r+1} (t)$ is $\{ t^{r-k}\tilde{\lie{b}}_i| \lie{b}_i \in B_i , k\leq i\leq r\} \dot \cup \{ t^i\tilde{\lie{b}}_j | \lie{b}_j \in B_j, j\geq r+1, 0\leq i \leq r \}$.}
 \item The Lie algebra $\mf{h}^{r+1}(0)$ is precisely $\wmf{g}_{r+1}$ and is the Lie algebra of $H^{r+1}(0)$. 
    \end{enumerate}
\end{prop}

\noindent
{\bf Proof}: (1) follows Prop. \ref{prop:jet} (2). Let us now prove (2). 

Start with the basis $B$ constructed above. From Proposition~\ref{prop:hatbasis} and Lemma~\ref{lemma:gammahat}, $\{\wmf{\lie{b}} | \lie{b} \in B\}$ are a basis of $\wmf{g}_0 = \mf{g}$. Recall that for $\lie{b}(t)=\gamma (t)\lie{b} \gamma(t)^{-1}$, the leading term is $t^{d(\lie{b})}\widehat{\lie{b}}$. So for each $\lie{b} \in B$, let us construct the normalized operator
$$ \tilde{\lie{b}}(t) = t^{-d(\lie{b})}\lie{b}(t) \in \mf{g}\otimes {\cal O}$$
and let $\tilde{B}= \{\tilde{\lie{b}}| \lie{b}\in B \}$. Then $\tilde{B}$ forms an ${\cal O}$-basis of $\mf{g} \otimes {\cal O}$. 

Furthermore note that by Proposition~\ref{prop:gry}(2), 
the leading terms $\widehat{\lie{b}y}$ form a basis of $\widehat{T_y}$ and are linearly independent. We then define:
$$ W_{\lie{b}}(t) : = t^{-d(\lie{b}y)}\gamma(t) \lie{b}y, \: \: \lie{b} \in B',$$

Note that these are elements in $V\otimes {\cal O}$ with leading term $\widehat{\lie{b}y}$. Thus $\tilde{\lie{b}}(t)y(t) = t^{\delta(\lie{b})}W_{\lie{b}}(t)$. 

Now, given $\lie{u}(t) \in \mf{h}^r (t)$, write
\[ \lie{u}(t) = \sum_{\lie{b} \in B} f_{\lie{b}}(t) \tilde{\lie{b}}(t), f_{\lie{b}}(t) \in {\cal O}. \]

 Hence,
\[ \lie{u}(t)y(t) = \sum_{\lie{b} \in B'} f_{\lie{b}}(t)t^{\delta(\lie{b})}W_{\lie{b}}(t).\]
So if $\lie{u}(t) \in \mf{h}^{r}(t)$, then $f_{\lie{b}}(t)t^{\delta(\lie{b})} \in t^{r+1}{\cal O}.$
So $f_{\lie{b}}(t) \in t^{\operatorname{max}(0, r+1-\delta(\lie{b}))}$. Thus (2) is proved. Now for (3), evaluating at zero gives us an element in the span of $\{ \hat{\lie{b}} | \delta(\lie{b}) \geq r+1\}$. 
This is precisely $\wmf{g}_{r+1}$. Conversely, given any $\lie{g}\in \mf{g}_r$ such that $\wl{g}\in \wmf{g}_r$, the expansion $\lie{g}(t)= t^{d(\lie{g})}g_d + \ldots + t^{d(\lie{g})+r}g_{d+r}+\ldots$ provides us the element $\lie{g}'(t)=g_0 +t^1 g_1 +\ldots + t^r g_r \in \mf{g}\otimes \C[r]$ such that $\lie{g}'(t)y^r =0$. This proves the proposition. $\Box$

\begin{corollary} \label{lemma:stabilizer}
The cancellation alegbra $\wmf{g}_2$ equals $\mf{h}_{\overline{y_1}}$. 
\end{corollary}

\noindent
{\bf Proof}: Now it is clear that $\wmf{g}_2  \subseteq \mf{h}_{\overline{y_1}}$. For the converse, let $\lie{g}_0\in \mf{h}_{\overline{y_1}}$. Then there is a $\lie{g}_1$ such that $\lie{g}_0 y_1 +\lie{g}_1 y_0=0$. We already have $\lie{g}_0 y_0=0$. Thus $\lie{g}(t)=\lie{g}_0 +t\lie{g}_1 \in \mf{h}^2 (t)$. But my Prop. \ref{prop:lie}(3), $\lie{g}_0 = \lie{g}(0)\in \wmf{g}_2$. This proves the corollary. $\Box$

\paragraph{Remark.}
The results above construct an algebraic subgroup of $H=G_z$
whose Lie algebra is $\wmf{k}$ by evaluating a sufficiently
high-order jet stabilizer $H^r$ at zero. We leave the reader with the following question : Do there exist algebraic groups
$A\subseteq H^r$, for some $r$, such that (i) $
Lie(\pi^0 (A))=\wmf{k}
$
but (ii) $\pi^0 (A)$ ranges over multiple components of $H$?
\eat{\paragraph{Question.}
Let $\mathcal O=\mathbb C[[t]]$ and let
$y(t)=\gamma(t)y\in V\otimes_{\mathbb C}\mathcal O$, with $y(0)=z$.
Let $\wmf{k}$ be the limiting stabilizer Lie algebra constructed above.

Do there exist a connected finite-dimensional algebraic groups
$A$ over $\mathbb C$ and an injective homomorphism
\[
\iota:A\longrightarrow G(\mathcal O)
\]
such that every element of $\iota(A)$ stabilizes $y(t)$ and
\[
\Lie\!\left(\operatorname{ev}_0(\iota(A))\right)=\wmf{k}?
\]
Here algebraicity means that, for every $r\geq 0$, the induced map
\[
A\longrightarrow G\bigl(\mathbb C[t]/(t^{r+1})\bigr)
\]
is a morphism of algebraic groups.
}

\section{Examples}
\begin{ex}
    Let $V=\C^{3\times 3}$ and $G=GL_3 (\C)$ with the action as left multiplication. In other words, for $g\in G$ and $x\in V$, we have $\rho (g)(x)=gx$. Let $\lambda (t)$ and $y$ be as given below.
\[ \lambda(t)= \left[ \begin{array}{ccc}
1 & 0 & 0 \\
0 & t & 0 \\
0 & 0 & t^2 \\
\end{array} \right] 
\: 
\: 
y=\left[ \begin{array}{ccc}
1 & 0 & 0 \\
0 & 1 & 0 \\
0 & 0 & 1 \\
\end{array} \right] 
\: \: 
z=\left[ \begin{array}{ccc}
1 & 0 & 0 \\
0 & 0 & 0 \\
0 & 0 & 0 \\
\end{array} \right] 
\: \: 
\mf{h}=\left[ \begin{array}{ccc}
0 & * & * \\
0 & * & * \\
0 & * & * \\
\end{array} \right]
\]
We see that $y\stackrel{\lambda}{\rightarrow}z$. Let $(e_{ij})$ be the standard basis for $V$, then we have:
\[ y(t)=\lambda(t)y = e_{11}+te_{22} +t^2 e_{33} \]
We also see that $\mf{g}=End(\C^3)$ consists of $3\times3$-matrices. The stabilizer $\mf{k}=\mf{g}_y =(0)$ while $dim(\mf{h})=6$ and consists of elements $\lie{g}\in \mf{g}$ with non-zero entries in columns $2$ and $3$ (see above).

Let a typical element $\lie{g} \in \mf{g}$ be given by the $3\times 3$-matrix $(g_{ij})$. We then have $\lie{g}(t)$ and $\lie{gy}(t)$ as given below:
\[ \lie{g}(t)=
\left[ \begin{array}{ccc}
g_{11} & t^{-1}g_{12} & t^{-2}g_{13}\\
t^1 g_{21} & g_{22} & t^{-1}g_{23}\\
t^2 g_{31} & t^1 g_{32} & g_{33}\\
\end{array} \right]
\: \:
\lie{gy}(t)=
\left[ \begin{array}{ccc}
g_{11} & g_{12} & g_{13}\\
t^1 g_{21} & t^1 g_{22} & t^{1}g_{23}\\
t^2 g_{31} & t^2 g_{32} & t^2g_{33}\\
\end{array} \right]
\]
We see that $-2\leq d(\lie{g})\leq 2$, while $0\leq d(\lie{g}y)\leq 2$. Thus, $\mf{g}^{-2}=\mf{g}^{(0)}=\mf{g}$. 
\[ \mf{g}^{-1}=
\left[ \begin{array}{ccc}
g_{11} & g_{12} & 0\\
g_{21} & g_{22} & g_{23}\\
 g_{31} & g_{32} & g_{33}\\
\end{array} \right]
\:
\: 
\mf{g}^{0}=
\left[ \begin{array}{ccc}
g_{11} & 0 & 0\\
g_{21} & g_{22} & 0\\
 g_{31} & g_{32} & g_{33}\\
\end{array} \right]
\]
The others - $\mf{g}^1$ and $\mf{g}^2$ - are as expected. 
\[ 
\mf{g}^{(1)}=
\left[ \begin{array}{ccc}
0 & 0 & 0\\
g_{21} & g_{22} & g_{23}\\
 g_{31} & g_{32} & g_{33}\\
\end{array} \right]
\mf{g}^{(2)}=
\left[ \begin{array}{ccc}
0 & 0 & 0\\
0 & 0 & 0\\
 g_{31} & g_{32} & g_{33}\\
\end{array} \right]
\]
Let $d_{ij}=dim(\mf{g}^{ij})$ and $e_{ij}=|B_{ij}|=dim(\mf{g}^{i,j}/(\mf{g}^{i+1,j}+\mf{g}^{i,j+1}))$. These tables are given below:
\[
(d_{ij})=\begin{array}{|r||c|c|c|c|}\hline 
i\downarrow j\rightarrow& 0 & 1 & 2 & 3\\ \hline \hline 
-2 &9 &6 & 3 &0\\ \hline
-1 & 8& 6& 3 &0\\ \hline
0 & 6&5 &3 &0\\ \hline
1 & 3& 3& 2&0\\ \hline
2 & 1&1 &1 &0\\ \hline
3 & 0 & 0 & 0 & 0 \\ \hline 
\end{array} \: \: 
(e_{ij})=\begin{array}{|r||c|c|c|}\hline 
i\downarrow j\rightarrow& 0 & 1 & 2 \\ \hline \hline 
-2 &1 &0 & 0\\ \hline
-1 & 1& 1& 0\\ \hline
0 & 1&1 &1 \\ \hline
1 & 0& 1& 1\\ \hline
2 & 0&0 &1 \\ \hline
\end{array}
\]
We may compute $B$ from the above. A summary of the basis $B$ appears below:
\[ 
\left[ \begin{array}{ccc}
b_{0,0} & b_{-1,0} & b_{-2,0}\\
b_{1,1} & b_{0,1} & b_{-1,1}\\
 b_{2,2} & b_{1,2} & b_{0,2}\\ \end{array} \right] \]
 Thus, as expected, $dim(\wmf{g}_1)=6$ and $\wmf{g}_1 = \mf{h}$, while $\wmf{g}_2 = \mf{h}_{\overline{y_1}}$. 
We may also compute $T_r=(\widehat{Ty}_r)$ as the spaces spanned by $\wl{by}$, for $\lie{b} \in B_{i,j}$, for $j-i=r$. These spaces are summarized below:
\[ 
\left[ \begin{array}{ccc}
T_0 & T_1 & T_2 \\
T_0 & T_1 & T_2 \\
T_0 & T_1 & T_2 \\
\end{array} \right] \]
Let $w=y_0 +y_1 = e_{11}+e_{22}$. Then, we have $\overline{O(z)}\subsetneq \overline{O(w)}\subsetneq \overline{O(y)}$. Note that $T_0$ is the tangent space $T_z$, while $T_1$ is the space of normals for the containment $\overline{O(z)}\subsetneq \overline{O(w)}$ and $T_2$ is the space of normals for the containment $\overline{O(w)}\subsetneq \overline{O(y)}$.
\end{ex}
\begin{ex}
    Let $V, G, y$ be as above. However, we consider:
\[ \mu(t)= \left[ \begin{array}{ccc}
1 & 0 & 0 \\
0 & t & 0 \\
0 & 0 & t \\
\end{array} \right] 
\: 
\: 
y=\left[ \begin{array}{ccc}
1 & 0 & 0 \\
0 & 1 & 0 \\
0 & 0 & 1 \\
\end{array} \right] 
\: \: 
z=\left[ \begin{array}{ccc}
1 & 0 & 0 \\
0 & 0 & 0 \\
0 & 0 & 0 \\
\end{array} \right] \]
We then have $y \stackrel{\mu}{\rightarrow} z$ as before, but $y(t)$ is now different:
\[ y(t)=\mu(t)y = e_{11}+t(e_{22} + e_{33}) \]   
We may compute $B^{\mu}$ from the above. A summary of the basis $B$ appears below:
\[ 
\left[ \begin{array}{ccc}
b_{0,0} & b_{-1,0} & b_{-1,0}\\
b_{1,1} & b_{0,1} & b_{0,1}\\
 b_{1,1} & b_{0,1} & b_{0,1}\\ \end{array} \right] \]
 Thus, as before, $dim(\wmf{g}^{\mu}_1)=6$ and $\wmf{g}^{\mu}_1 $ equals $\mf{h}$. But $\wmf{g}^{\mu}_2 =(0)$.
We may also compute $T^{\mu}_r=(\widehat{Ty}^{\mu}_r)$ as before. These are summarized below:
\[ 
\left[ \begin{array}{ccc}
T^{\mu}_0 & T^{\mu}_1 & T^{\mu}_1 \\
T^{\mu}_0 & T^{\mu}_1 & T^{\mu}_1 \\
T^{\mu}_0 & T^{\mu}_1 & T^{\mu}_1 \\
\end{array} \right] \]
Thus $\mu$ does {\bf not visit} any intermediate orbit and $T^{\mu}_1$ is the normal space to the containment $\overline{O(z)}\subsetneq \overline{O(y)}$. 
\end{ex}
\eat{\begin{ex}
    Let $V=\C^{2\times 2}$ and $G=GL_2 (\C)$ acting on $V$ be left multiplication. Let $\lambda (t)$ and $y$ be as given below.
\[ \lambda(t)= \left[ \begin{array}{cc}
1 & 0  \\
0 & t  \\
\end{array} \right] 
\: 
\: 
y=\left[ \begin{array}{cc}
1 & 0  \\
1 & 1  \\
\end{array} \right] 
\: \: 
z= \left[ \begin{array}{cc}
1 & 0  \\
0 & 0  \\
\end{array} \right]
y(t)=\left[ \begin{array}{cc}
1 & 0  \\
t & t  \\
\end{array} \right]\]
\[ \lie{g}(t)=
\left[ \begin{array}{cc}
g_{11} & t^{-1}g_{12} \\
t^1 g_{21} & g_{22} \\
\end{array} \right]
\: \:
\lie{g}(t)y(t)=\lie{gy}(t)=
\left[ \begin{array}{cc}
(g_{11}+g_{12}) & g_{12} \\
t^1 (g_{21}+g_{22}) & t^1 g_{22} \\
\end{array} \right]
\]

\end{ex}}
\begin{ex}
    Let $V=Sym^2 (\C^2)$ acted upon by $G=GL_2 (\C)$. Consider the point $p=x_1^2 +x_2^2$ and $q=x_1^2$ obtained by $\lambda (t)$ given below. We have $\mf{k}$ and $p(t)$ as shown below.
\[ \lambda(t)= \left[ \begin{array}{cc}
1 & 0  \\
0 & t  \\
\end{array} \right] 
\: 
\: 
\kappa = \C \left[ \begin{array}{cc}
0 & 1  \\
-1 & 0  \\
\end{array} \right]
\: \: p(t)=x_1^2 +t^2 x_2^2 \: \: 
\kappa (t) = \left[ \begin{array}{cc}
0 & t^{-1}  \\
-t & 0  \\
\end{array} \right]
\]
For a typical element $\lie{g}\in \mf{g}$, we have $\lie{g}(t)$ and $\lie{gp}(t)$ as below:
\[ \lie{g}(t) \left[ \begin{array}{cc}
g_{11} & t^{-1} g_{12}  \\
tg_{21} & g_{22}  \\
\end{array} \right]
\: \: \lie{gp}(t)=2g_{11}x_1^2 + 2(g_{12}+g_{21})tx_1 x_2   +2g_{22} t^2 x_2^2 \]
As a result, we have:
\[ \mf{g}^{-1}=\mf{g}, \: \: \mf{g}^{0} = \left[ \begin{array}{cc}
g_{11} &  0  \\
g_{21} & g_{22}  \\
\end{array} \right], \: \:  
\mf{g}^{1} = \left[ \begin{array}{cc}
0 &  0  \\
g_{21} & 0  \\
\end{array} \right]
\]
\[ \mf{g}^{(0)}=\mf{g}, \: \: \mf{g}^{(1)} = \left[ \begin{array}{cc}
0 &  g_{12}  \\
g_{21} & g_{22}  \\
\end{array} \right], \: \:  
\mf{g}^{(2)} = \left[ \begin{array}{cc}
0 &  g_{12}  \\
-g_{12} & g_{22}  \\
\end{array} \right]
\]
Let $d_{ij}=dim(\mf{g}^{ij})$ and $e_{ij}=|B_{ij}|=dim(\mf{g}^{i,j}/(\mf{g}^{i+1,j}+\mf{g}^{i,j+1}))$. These tables are given below:
\[
(d_{ij})=\begin{array}{|r||c|c|c|c|}\hline 
i\downarrow j\rightarrow& 0 & 1 & 2 & 3\\ \hline \hline 
-1 & 4& 3& 2 &1\\ \hline
0 & 3&2 &1 &0\\ \hline
1 & 1& 1& 0&0\\ \hline
2 & 0& 0& 0 &0\\ \hline 
\end{array} \: \: 
(e_{ij})=\begin{array}{|r||c|c|c|}\hline 
i\downarrow j\rightarrow& 0 & 1 & 2 \\ \hline \hline 
-1 & 0& 0& 0\\ \hline
0 & 1&0 &1 \\ \hline
1 & 0& 1& 0\\ \hline
\end{array}
\]
This gives us the basis $B$:
\[ \lie{b}_{0,0}=\left[ \begin{array}{rr}
1 & 0  \\
0 & 0  \\
\end{array} \right], \: \:
\lie{b}_{1,1}=\left[ \begin{array}{rr}
0 & 0  \\
1 & 0  \\
\end{array} \right], \: \:
\lie{b}_{0,2}=\left[ \begin{array}{rr}
0 & 0  \\
0 & 1  \\
\end{array} \right], \: \:
\lie{b}_{-1,3}=\left[ \begin{array}{rr}
0 & 1  \\
-1 & 0  \\
\end{array} \right]
\]
I think the element KV meant was:
\[  t\lie{b}_{-1,3}(t)=t\lie{k}(t)= \left[ \begin{array}{cc}
0 & 1  \\
-t^2 & 0  \\
\end{array} \right]
\]
I think the exponential of this element, which should be in $H_2 ({\cal O})$ is:
\[   h(t)=\left[ \begin{array}{cc}
\cos t & \sin t/t \\
-t\sin t & \cos t  \\
\end{array} \right] \]
One may check:
\[ (x_1 \cos t-t x_2 \sin t)^2 + t^2 ( x_1 t^{-1} \sin t + x_2 \cos t)^2 = p(t)\]
Note that:
\[   h(0)=\left[ \begin{array}{cc}
1 & 1 \\
0 & 1
\end{array} \right] \]
\end{ex}

\begin{ex}[Binary cubics: stabilizers, paths, and cancellation]
\label{ex:binary-cubics-combined}
\label{ex:binary-cubics-representation}
Let $G=GL_2(\C)$ and $V=\operatorname{Sym}^3(\C^2)$ under the standard action $\rho :G \rightarrow GL(V)$. We use the ordered basis
\[
\mathcal E=(x^3,x^2y,xy^2,y^3).
\]
Let $p = x^2y, q = x^3$. In this example we compute the filtration algebras for $x^2y \stackrel{\mu}{\rightarrow} x^3$, under different 1-PF's 
$\mu$.  

We use column vectors for coefficients in $\mathcal E$.
The action of the Lie algebra $\mathfrak{gl}_2$ on a polynomial P is given by:
$$
\lie{g}=\begin{pmatrix}g_{11}&g_{12}\\g_{21}&g_{22}\end{pmatrix}
\in\mf{g}=\mathfrak{gl}_2,
\qquad
\lie{g}\cdot P
=(g_{11}x+g_{21}y)\partial_xP
+(g_{12}x+g_{22}y)\partial_yP.
$$
The matrix of the representation is 
\[\rho_*(\lie{g})=
\begin{pmatrix}
3g_{11}&g_{12}&0&0\\
3g_{21}&2g_{11}+g_{22}&2g_{12}&0\\
0&2g_{21}&g_{11}+2g_{22}&3g_{12}\\
0&0&g_{21}&3g_{22}
\end{pmatrix}.
\]
Let $\gamma (t)$ be as below. Then we have:
\[ p(t)=\gamma (t)p= (x^2)(ty+x)=x^3 +t x^2 y\] Thus $p_0=x^3$ and $p\stackrel{\gamma}{\rightarrow} q$. Note that $p_1 = x^y \in T_q$, thus the path approaches tangential to the orbit $O(q)$. 
\[
\gamma(t)=\begin{pmatrix}1&1\\0&t\end{pmatrix},\qquad
R(t)=\rho(\gamma(t))=
\begin{pmatrix}
1&1&1&1\\
0&t&2t&3t\\
0&0&t^2&3t^2\\
0&0&0&t^3
\end{pmatrix}.
\]
Set 
\[
L=-g_{11}+g_{12}-g_{21}+g_{22},\qquad
\ell_0=2g_{11}+g_{12}+2g_{21}+g_{22},\qquad
\ell_1=2g_{11}+4g_{21}+g_{22}.
\]
We may compute: 
$$
\lie{g}(t)=
\begin{pmatrix}
g_{11}+g_{21}&L/t\\
t g_{21}&g_{22}-g_{21}
\end{pmatrix},
\qquad
\gamma(t)(\lie{g}p)
=\ell_0x^3+t\ell_1x^2y+2t^2g_{21}xy^2.
$$ 
The filtrations are given by
$$
\mf{g}^i=
\begin{cases}
\mf{g},&i\le-1,\\
\{\lie{g}:L=0\},&i=0,\\
\{\lie{g}:L= g_{11}+g_{21}= g_{22}-g_{21}=0\},&i=1,\\
0,&i\ge2
\end{cases}
\quad
\mf{g}^{(j)}=
\begin{cases}
\mf{g},&j\le0,\\
\{\lie{g}:\ell_0=0\},&j=1,\\
\{\lie{g}:\ell_0=\ell_1=0\},&j=2,\\
\{\lie{g}:\ell_0=\ell_1=g_{21}=0\}
=\C\lie{\kappa},&j\ge3.
\end{cases}
$$
Here $\lie{\kappa} = E_{11} - 2E_{22}$ in standard matrix coordinates.

As before, define
$d_{ij}=\dim(\mf{g}^{i,j})$,
and
$e_{ij}=|B_{ij}|=
\dim\left(
\mf{g}^{i,j}/
({\mf{g}^{i+1,j}+\mf{g}^{i,j+1})}
\right)
$. The tables below record these dimensions.
The column $j=3$ records the stabilizer $\kappa$ of $p$.
\[
(d_{ij})=
\begin{array}{|r||c|c|c|c|}
\hline
i\downarrow\,j\rightarrow & 0 & 1 & 2 & 3 \\
\hline\hline
-1 & 4 & 3 & 2 & 1 \\ \hline
 0 & 3 & 2 & 1 & 0 \\ \hline
 1 & 1 & 1 & 0 & 0 \\ \hline
 2 & 0 & 0 & 0 & 0 \\ \hline
\end{array}
\qquad
(e_{ij})=
\begin{array}{|r||c|c|c|c|}
\hline
i\downarrow\,j\rightarrow & 0 & 1 & 2 & 3 \\
\hline\hline
-1 & 0 & 0 & 0 & 1 \\ \hline
 0 & 1 & 0 & 1 & 0 \\ \hline
 1 & 0 & 1 & 0 & 0 \\ \hline
\end{array}
\]
Consider the following Lie algebra elements:
\[
\lie{v}_{00}=\begin{pmatrix}1&1\\0&0\end{pmatrix},\quad
\lie{v}_{11}=\begin{pmatrix}-1&-1\\1&1\end{pmatrix}, \quad
\lie{v}_{02}=\begin{pmatrix}-1&2\\1&-2\end{pmatrix},\quad
\lie{\kappa}=\lie{v}_{-1,3}=\begin{pmatrix}1&0\\0&-2\end{pmatrix}
\]
We have,
\[
\lie{v}_{00}(t)=\begin{pmatrix}1&0\\0&0\end{pmatrix},\quad
\lie{v}_{11}(t)=\begin{pmatrix}0&0\\t&0\end{pmatrix}, \: 
\lie{v}_{02}(t)=\begin{pmatrix}0&0\\t&-3\end{pmatrix},\quad
\lie{\kappa}(t)=\begin{pmatrix}1&-3/t\\0&-2\end{pmatrix}\quad
\]
The four leading terms are independent. They are $\lie{E}_{11}$,
$\lie{E}_{21}$, $\lie{E}_{22}$, and $-3\lie{E}_{12}$, respectively. Moreover: 
$$
\gamma(t)(\lie{v}_{00}p)=3x^3+2t x^2y,\quad
\gamma(t)(\lie{v}_{11}p)=3t x^2y+2t^2xy^2, \quad
\gamma(t)(\lie{v}_{02}p)=2t^2xy^2,\quad and \ 
\gamma(t)(\lie{\kappa}p)=0
$$

Thus, we have:
\[
B_{0,0}=\{\lie{v}_{00}\},\quad
B_{1,1}=\{\lie{v}_{11}\}, \quad
B_{0,2}=\{\lie{v}_{02}\}, \quad
B_{-1,3}=\{\lie{\kappa}\},\quad
\]
Note that $dim(T_p)=3$ while $dim(T_q)=2$. We have $\widehat{(Tp)}_0=\C \{ x^3, x^2y\}$, $\widehat{(Tp)}_1=0$ and $\widehat{(Tp)}_0=\C \{ xy^2\}$. 
These witnesses the inclusion $\overline{O(q)}\subsetneq \overline{O(q)}$. Even though $p(t)$ approaches $q$ tangential to the orbit, $\widehat{(Ty)}_2=\C \{xy^2\}$ is normal to $T_q$.
We also have:
\[
\mf{g}=\wmf{g}_0 \supset
\wmf{g}_1=\wmf{g}_2
=\langle\lie{E}_{12},\lie{E}_{22}\rangle=\mf{h}\supset 
\wmf{g}_3=\C\lie{E}_{12}=\wmf{k}
\]
Now, let us construct the group element. We have $\kappa'(t)$ and $e^{\kappa'(t)}$ as follows:
\[ 
\kappa ' (t)=t\lie{\kappa}(t)=\begin{pmatrix}t&-3\\0&-2t\end{pmatrix}\quad
g(t)=e^{\kappa ' (t)}=\begin{pmatrix}e^t&(e^{-2t}-e^{t})/t\\0&e^{-2t}\end{pmatrix}
\]
One may check that $g(t)\cdot p(t)=p(t)+O(t^2)$.  
\[ \begin{array}{rcl}
g(t) p(t) &=& e^{3t} x^3 + e^{2t}x^2 (te^{-2t}y + (e^{-2t} -e^t)x \\
&=& e^{3t} x^3 + tx^2 y +(1-e^{3t})x^3\\
&=& x^3 +tx^2 y
\end{array}\]
Note that:
\[ 
g(0)= \begin{pmatrix} 1 & -3 \\ 0 & 1 \end{pmatrix} \]
which stabilizes $p(0)=x^3$.
\end{ex}

\eat{\textcolor{green}{ Next example can go ? It is concerned with $\gamma^{-1}(t) \dot{\gamma(t)}$ and this paper does not really talk about this? 
\begin{ex} \label{ex:HM}
This is Example 7. Let $G=SL_2 (\C)$ and $V=Sym^3(\C)$. Let $p=x^2 y$ and $q=x^3$. Note that $G_p$ is trivial and $G_q$ is the unipotent group of dimension $1$. Let $\gamma (t)$ be as given below:
\[ \gamma (t)\left[ \begin{array}{c}
x\\
y \end{array} \right]
=\left[ \begin{array}{cc}
t & 0 \\
t^{-2} & t^{-1} \end{array} \right]\left[ \begin{array}{c}
x\\
y \end{array} \right]\]
We have:
\[ p([x,y]^T)(t)=p(\gamma(t)[x, y]^T) = (tx)^2 (t^{-2}x +t^{-1}y)=x^3 +tx^2 y \]
Thus we have $p \stackrel{\gamma}{\rightarrow} q$ with $e=1$ and $x^2 y$ is the tangent vector, which is in $T_q O(q)$. 
The provisional construction $\mf{v}(t)$ is given below:
\[ \mf{v}(t)=\dot{\gamma}(t) \gamma^{-1}(t)= \left[ \begin{array}{cc}
1 & 0 \\
-2t^{-3} & -t^{-2} \end{array} \right]
\left[ \begin{array}{cc}
t^{-1} & 0 \\
-t^{-2} & t \end{array} \right]
= 
\left[ \begin{array}{cc}
t^{-1} & 0 \\
-t^{-4} & -t^{-1} \end{array} \right]
\]
On the other hand, we have:
\[ \mf{v}'(t)=\gamma^{-1}(t) \dot{\gamma(t)}= 
\left[ \begin{array}{cc}
t^{-1} & 0 \\
-t^{-2} & t \end{array} \right]
\left[ \begin{array}{cc}
1 & 0 \\
-2t^{-3} & -t^{-2} \end{array} \right]
= 
\left[ \begin{array}{cc}
t^{-1} & 0 \\
-3t^{-2} & -t^{-1} \end{array} \right] \in \mf{g}\otimes k
\]
\textcolor{red}{for some dual reasons the role of $\mf{v}$ and $\mf{v}'$ are interchanged.} This $\mf{v}(t)p(t)=\dot{p}(t)=x^2y$. 
We then have:
\[ \begin{array}{rcl} 
\mf{v}(t)p &=& t^{-1}x \frac{\partial p}{\partial x}
-3t^{-2} x \frac{\partial p}{\partial y} - t^{-1} y \frac{\partial p}{\partial y} \\
&=& -t^{-4} x^3 + t^{-1} x^2 y\\ 
\end{array}\]
Thus $\widehat{\mf{v}'(t)p}= x^3$. Then $\gamma (t)(x^3)=t^3 x^3$. However the tuple $(p, \mf{v}p)\stackrel{\gamma}{\rightarrow} (q,p)$, just that $\mf{v}p\in V\otimes k$! 
\end{ex}
}

\textcolor{blue}{Removed non-ample example} 

\textcolor{red}{The next portions are about the normal cone - can be removed}

***** From here *******
The orbit closure of $y$ is given by the vanishing $\Delta$, of the discriminant of $a x^3 + b x^2 y + c xy^2 + dy^3$,
\[ \Delta
=
b^2c^2-4ac^3-4b^3d-27a^2d^2+18abcd. \]
The orbit closure of $O(z)$ is the variety
$Z=\overline{O(x^3)}$, the cone of perfect cubes. On the open set $a\neq0$, we can write a perfect cube as
\[
A(x+s y)^3.
\]
Its coefficients are
\[
a=A,\qquad
b=3As,\qquad
c=3As^2,\qquad
d=As^3.
\]
Since
\[
s=\frac{b}{3a},
\]
we obtain
\[
c=\frac{b^2}{3a},
\qquad
d=\frac{b^3}{27a^2}.
\]
Introducing local normal coordinates, $C=c-\frac{b^2}{3a},
D=d-\frac{b^3}{27a^2}$, locally along $Z$ we have $I_Z=(C,D).$ In terms of these local coordniates we have 
\[
\Delta
=
-4aC^3
-3b^2C^2
+18abCD
-27a^2D^2.
\]
So the initial ideal of $I_Y$ in $I_Z$ is given $-3b^2C^2+18abCD-27a^2D^2
=-3(bC-3aD)^2$. The coordinate ring of the normal cone to $Z$ is the 
the associated graded ring $\operatorname{gr}_{I_Z}\C[V]$. On the other hand
$(bC-3aD)^2$ cuts out the normal cone $C_{Z}Y$. At $z=x^3$ this becomes $D^2=0$, and so $(C_ZY)_z = \Spec\frac{\C[C,D]}{(D^2)}.$
Its reduced support is $D=0$, i.e the $C$ axis, namely the line $\C x^2y \subset N_z$.  So the exceptional fiber of 
$\operatorname{Bl}_ZY\longrightarrow Y$ over $z$ has reduced support consisting of the single point $\overline{xy^2}$. 

****   TILL HERE     ******

\textcolor{blue}{Do we need one more example, in the null cone? Not needed I think. We should 
remember this when we want to understand the blow-up.}
\eat{
********* From HERE *******

\begin{ex}
{\bf Example 3 (ample)}: $\gamma(t)=\begin{psmallmatrix}1&0\\1&t\end{psmallmatrix}$

Row action; conjugation $\lie{g}(t)=\gamma(t)\lie{g}\gamma(t)^{-1}$.
$y(t)=tx^2y+2txy^2+ty^3$ (no constant term), so $z=0$.
\[
\mf{h}=\operatorname{Stab}(0)=\mf{gl}_2\quad\text{(trivially),}\qquad\dim=4.
\]

\noindent\textit{Ampleness.} $w(t)$ collected by power of $t$:
\[
t^0:b(x+y)^3,\qquad t^1:a(x^2y+2xy^2+y^3),\qquad t^2:c(xy^2+y^3).
\]
$(a,b,c)=(0,1,0),(1,0,0),(0,0,1)$ give leading terms
$(1,3,3,1),(0,1,2,1),(0,0,1,1)$ in the $(x^3,x^2y,xy^2,y^3)$-basis --
rank $3$, matching $\dim T_yO(y)$. \textbf{Ample.}

\noindent\textit{$\mf{g}_1$ and its leading terms.}
$\lie{g}(t)=\begin{psmallmatrix}c_1-c_2/t & c_2/t\\ c_1-c_2/t+c_3t-c_4 & c_2/t+c_4\end{psmallmatrix}$.
$\delta(\lie{g})-d(\lie{g})\ge1$ on every stratum outside $\mf{k}$: $c_2\ne0$
($d=-1,\delta=0$); $c_2=0$, $2c_1+c_4\ne0$ ($d=0,\delta=1$); $c_2=0,2c_1+c_4=0$,
$c_3\ne0$ (deeper, $\delta-d=2$); with $\mf{k}=\operatorname{span}(E_{11}-2E_{22})$
itself giving $\delta=\infty$. So
\[
\mf{g}_1=\mf{gl}_2\quad\text{(every nonzero element satisfies }\delta-d\ge1\text{)},\qquad
\dim(\operatorname{span}\mf{g}_1)=4.
\]
Leading terms across these strata span all of $\mf{gl}_2$:
\[
\widehat{\mf{g}_1}=\mf{gl}_2=\mf{h},\qquad\dim=4.
\]
$\widehat{\mf{g}_1}=\mf{h}$ exactly, \emph{and}
$\dim\widehat{\mf{g}_1}=4=\dim(\operatorname{span}\mf{g}_1)$:

$Y$ is itself a cone, with vertex $0$, since $\Delta$ is homogeneous of
degree $4$. Blowing up $z=0$ in $V$ gives exceptional divisor
$E_0=\PP(V)\cong\PP^3$, and since $Y$ is already a cone through $0$,
\[
\widetilde Y\cap E_0 = \PP(Y) \subset \PP^3.
\]
The strict transform meets $E_0$ along the entire projectivization
of $Y$ itself. Since $\Delta$ is homogeneous of degree 4, no conditions are imposed in degrees 1,2,and 3, and so $Y^1(0)=V$,$Y^2(0)=V \times V$ and $Y^3(0) = V \times V \times V$. $Y^4(0)$'s imposes the condition $\Delta(v_1)=0$,
, i.e. \ $v_1\in Y$.

The three leading terms from the discussion at the start of this example are
\[
\hat w_1=(x+y)^3,\qquad \hat w_2=2y(x+y)^2,\qquad \hat w_3=2y^2(x+y).
\]
Now all of them have multiple roots, so all of them satisfy $\Delta(\hat w_i)=0$.
So all three meet $\widetilde Y\cap E_0=\PP(Y)$.
But they meet
at different depths of $\PP(Y)$'s singular locus.
$\hat w_1=(x+y)^3$ is a perfect cube, i.e.\
$[\hat w_1]\in\PP(Z)$, the deepest singular stratum of $\PP(Y)$.
$\hat w_2,\hat w_3$ have only \textbf{double} roots -- they land on
the smooth locus of $\PP(Y)$, away from $\PP(Z)$.
\end{ex}
}
\textcolor{red}{This is SYM4- again if ever we understand the blowup using the cancellation, then we can revive this }

****** FROM HERE    ******
\begin{ex}
$V=Sym^4(\C^2)$; $y=x^4+y^4$, $z=x^4$. Row action.

\[
T_yO(y)=\operatorname{span}(x^4,\,x^3y,\,xy^3,\,y^4),\qquad\dim=4.
\]
The $x^2y^2$ tangent is missing. Since $\gamma(t)$ is diagonal, this is ample.
\subsection*{$a$, $b$, $r_0$}
\[
a=-1,\qquad b=0,\qquad r_0=4.
\]

\subsection*{The bi-filtration table: $\dim\mf{g}^{i,j}$}
\[
\begin{array}{c|cccccc}
i\backslash j & 0 & 1 & 2 & 3 & 4 & 5\\ \hline
-1 & 4 & 3 & 2 & 2 & 1 & 0\\
0 & 3 & 2 & 1 & 1 & 1 & 0\\
1 & 1 & 1 & 0 & 0 & 0 & 0\\
2 & 0 & 0 & 0 & 0 & 0 & 0
\end{array}
\]
Columns $j=2$ and $j=3$ are \textbf{identical in every row}: no element of
$\mf{g}$ has $d(\lie{g}\cdot y)$ equal to exactly $2$.

\subsection*{The four basis elements}

All four are standard basis elements.

\bigskip\noindent\textbf{$E_{11}=\begin{psmallmatrix}1&0\\0&0\end{psmallmatrix}$},
\quad type $(0,0)$.
\[
E_{11}(t)=\begin{pmatrix}1&0\\0&0\end{pmatrix}=\widehat{E_{11}(t)}\quad(\text{constant in }t).
\]
\[
E_{11}\cdot y = 4x^4,\qquad \gamma(t)\cdot(E_{11}\cdot y) = 4x^4.
\]

\bigskip\noindent\textbf{$E_{21}=\begin{psmallmatrix}0&0\\1&0\end{psmallmatrix}$},
\quad type $(1,1)$.
\[
E_{21}(t)=\begin{pmatrix}0&0\\t&0\end{pmatrix},\qquad
\widehat{E_{21}(t)}=\begin{pmatrix}0&0\\1&0\end{pmatrix}=E_{21}.
\]
\[
E_{21}\cdot y = 4x^3y,\qquad \gamma(t)\cdot(E_{21}\cdot y) = 4t\,x^3y.
\]

\bigskip\noindent\textbf{$E_{12}=\begin{psmallmatrix}0&1\\0&0\end{psmallmatrix}$},
\quad type $(-1,3)$.
\[
E_{12}(t)=\begin{pmatrix}0&1/t\\0&0\end{pmatrix},\qquad
\widehat{E_{12}(t)}=\begin{pmatrix}0&1\\0&0\end{pmatrix}=E_{12}.
\]
\[
E_{12}\cdot y = 4xy^3,\qquad \gamma(t)\cdot(E_{12}\cdot y) = 4t^3\,xy^3.
\]

\bigskip\noindent\textbf{$E_{22}=\begin{psmallmatrix}0&0\\0&1\end{psmallmatrix}$},
\quad type $(0,4)$.
\[
E_{22}(t)=\begin{pmatrix}0&0\\0&1\end{pmatrix}=\widehat{E_{22}(t)}\quad(\text{constant in }t).
\]
\[
E_{22}\cdot y = 4y^4,\qquad \gamma(t)\cdot(E_{22}\cdot y) = 4t^4\,y^4.
\]

\bigskip
\noindent\textbf{Check.} $\{E_{11},E_{21},E_{12},E_{22}\}$ is the standard
basis, trivially rank $4$; $B\cap\mf{g}^{i,j}$ recovers exactly the table
above at every cell.

\textcolor{magenta}{To do: the Blow up of orbit closure of $x^4$ calculations.}
\end{ex}

**********    TILL HERE *******
}
\eat{
\section{Appendix: JoLT Review}

Let us review some of the notation and results of \cite{ASS25jolt}. 

\begin{defn}
 \label{def:Kt}
  Let $y,z,\gamma $ be as above. We define 
  \begin{enumerate}
      \item 
  $\widehat{\mf{k}}\subseteq {\mf{g}}$ as the $\C$-vector space generated by $\{ \widehat{\lie{k}} |\lie{k}\in {\mf{k}} \}$, the collection of leading terms of $\lie{k}(t)$, with $\lie{k} \in {\mf{k}}$.
  \item ${\mf{k}}(t)\subseteq {\mf{g}}'\otimes k$ as the $k$-vector space generated by all vectors $\lie{k} (t)$, for $\lie{k}\in {\mf{k}}$. 
  \item For an element $\lie{k}'\in \mf{g}\otimes k=\sum_{i\geq a} t^a \lie{g}_i$, with $\lie{g}_a \neq 0$, let $\widehat{\lie{k}'}=\mf{g}_a$. Let $\widehat{{\mf{k}}(t)}\subseteq {\mf{g}}$ as the $\C$-vector space generated by $\widehat{\lie{k}'}$, where $\lie{k}'\in {\mf{k}}(t)$. 
  
  \end{enumerate}
\end{defn}
It is shown that:
\begin{theorem}\label{thm:main}
 Let $y,z, \gamma $ be as in Equation~\ref{eqn1ps}.  Then $\widehat{\mf{k}}\subseteq \widehat{{\mf{k}}(t)}\subseteq {\mf{h}}$  and $\dim(\widehat{\mf{k}(t)})=\dim({\mf{k}})$. Moreover, $\widehat{{\mf{k}}(t)}\subseteq {\mf{h}}_{\overline{y_e}} $, the Lie algebra stabilizer within ${\mf{h}}$ of $\overline{y_e} \in \overline{N}$, where $y_e$ is the tangent of approach.  
\end{theorem}

\begin{defn}\label{Definition:assoc}
 Let $R =\oplus_{i\geq 0}  R_i$ be the associated graded ring for the ideal $I_z \subseteq \C [V]$. In other words, $R_i =I_z^i/ I_z^{i+1}$ (with $I_z^0 =\C [V]$). 
\end{defn}
Note that, for all $i$,  $R_i$ is a $G$-module.
For any ideal $I\subseteq I_z $, 
 set $I_i = I^i_z \cap I$ and $\overline{I}_i= (I^{i+1}_z +I_i)/I^{i+1}_z$ and let $\overline{I}=  \oplus_{i\geq 0} \overline{I}_i$. Then $\overline{I}$ is an ideal of $R$.
Note that for the $G$-invariant ideal $I_y \subseteq I_z$, the ideal $\overline{I_y}$ is a graded $G$-submodule of $R$. 

An easy consequence of the ring homomorphism $R_0\rightarrow R$ is the following:

\begin{prop} \label{proposition:projbundle}
    There is a morphism $Spec(R)\rightarrow \overline{O(z)}$. Thus, $X=\pi^{-1} 
    (O(z))$ is an open set of $Spec(R)$. Moreover, for any open set $U\subseteq O(z)$ we have the map $Spec(R_U)\rightarrow U$.  
\end{prop}

 Let $C$ be the affine variety $\overline{O(z)}\times V$ and $C^o$ be the open subset $O(z)\times V$ which is also the trivial bundle over $O(z)$. Let $TO(z)=\{(gz, gv)| g\in G, v\in T_z\}\subset C^o$ be the tangent bundle and $B^o=\{ (gz, v+gT_z)| g\in G, v\in V\}$, the quotient bundle. Note that $B^o \cong G\times^H \overline{N}$, the associated bundle over $O(z)$ with fiber $\overline{N}$.
 
 We denote a typical element of $B^o$ by $(w, \overline{v})$, where $w=gz\in O(z)$ and $\overline{v}=v+gT_z$. 
Let $\pi$ be the surjection $\pi: C^o \rightarrow B^o$ with fiber $gT_z$. For any open set $U\subseteq O(z)$, let $B^o_U =\{ (w,\overline{v}) \in B^o| w\in U\}$. 

\section{Proof of Commuting}
\textcolor{red}{To be written}
}
\eat{
\section{Tangents of approach}

\begin{defn}
\label{assume:nonesting} For $(z', \overline{v'})\in W$, we say that it satisfies the eigenvector condition if there is an $\lie{h}\in \mf{g}_{z'}$ such that $\lie{h} \overline{v'}=\alpha \overline{v'}$ with $\alpha \neq 0$. For $y\stackrel{\gamma}{\rightarrow} (z', \overline{v'})$, we say $\overline{v'}$ is a semisimple tangent of approach if $\ell_{\gamma} \overline{v'}=\alpha \overline{v'}$ with $\alpha \neq 0$. 
\end{defn}

\begin{theorem}
Let $\gamma(t)y = z+ y_e t^e+ \cdots +$ and suppose $\overline{y_e}$ is a semisimple tangent. Then there exists a 1-PS $\lambda $ and a $y' \in \overline{O(y)}$
 such that $y'\stackrel{\lambda}{\rightarrow}(z, \overline{y_e})$. \textcolor{red}{Should this not be $y' \in O(y)$ itself?} \end{theorem} 
 
 \begin{theorem}
    Let $G$ be reductive and $y \stackrel{\gamma}{\rightarrow} z$ and let the 1-PS $\lambda(t)$ stabilize $z$. Let $(-\lambda)$ denote that 1-PS $\lambda (t^{-1})$. Then either (A) there is a $y' \in (\overline{O(y)}-O(z))$ such that either (i) $\lambda$ stabilizes $y'$ or (ii) $y'\stackrel{\pm \lambda}{\longrightarrow} z$, or (B) $\lambda$ stabilizes a conjugate of $y$.
 \end{theorem}

Every element $x\in X$ is of the form (after conjugation by $G$) $z=(z, \overline{v})$, where $v\in T_z V$ and $\overline{v}$ is its image in $V\rightarrow \overline{N}$. For a tangent vector $x=(z,v)\in T_Z V$, let $J_{x}$ be the tangent ideal. Let $R[x]=R/J_x$ and $X_x =Spec(R[x])\cap X$. $R[x]$ is merely the coordinate ring of the $G$-orbit closure $\overline{O(x)}$.  

Let $R[y]=R/\overline{I}_y$ and $X_y =Spec(R[y])\cap X$. Note that if $x\in X_y$, for example, when $x=(z, \overline{y_e})$, the tangent of approach $J_x\supseteq \overline{I}_y$ and $X_x\subseteq X_y$. 

\begin{lemma} 
If $x$ satisfies the eigenvector condition, then $(z, \alpha \overline{v})\in O(x)$ for all $\alpha \in \C^*$. If $(\ell, \overline{v})$ is the pair as above then 
$[\ell, \mf{h}_{\overline{v}}]\subseteq \mf{h}_{\overline{v}}$.
\end{lemma}

We have the assumption:
\begin{assume}
\label{assume:ev}
We assume that for the limit $y\stackrel{\gamma}{\rightarrow}z$, the toa $\overline{y_e}$ satisfies the eigenvector condition.
\end{assume} 

\begin{prop}
Let $\gamma(t)y = z+y_e t^e+ \cdots +$. Suppose $\overline{y_e}$ satisfies the eigenvector condition with $\ell \cdot{y_e}= \alpha \overline{y_e}$, with $\alpha$ nonzero. Then there exists $\lambda = t^{\ell}$ a one-PS and $y' \in \overline{O(y)}$
 such that $y'\stackrel{\lambda}{\rightarrow}z$.    \textcolor{blue}{Anything stronger with $\gamma(t)$ satsfying the path condition?}
\end{prop}
\begin{lemma}
    \label{lemma:Iy}
 For $y$ as above, $dim(I_y)=dim(\overline{I}_y)$. Let $X_y =V(\overline{I}_y)\cap X$. Since $\overline{I}_y$ is homogeneous, if $(w,\overline{v})\in X_y$, then so is $(w, \alpha \overline{v})$ for any $\alpha \in \C^*$, and by closure, for $\alpha =0$ as well.    
\end{lemma}

The obvious question is: {\it Which $(z,\overline{v})$ are members of $F$?}. We have the first proposition (\cite{ASS25jolt}):

\begin{prop}
    \label{prop:jolt1}
    Let $\gamma (t)\subseteq G$ be a 1-parameter family such that:
   \[ \gamma (t) y = z + t^e v +\ldots  \] 
then $\overline{I}_z \supseteq J_{(z,\overline{v})} \supseteq \overline{I}_y$. Thus $(z,\overline{v})\in F$.  Moreover, if $\overline{v}\neq 0$, then we have an $\ell \in \mf{h}-\widehat{\mf{k}}$, arising from $\gamma$, such that $\overline{v}$ is an eigenvector of $\ell$. Thus $[\ell , \mf{h}_{\overline{v}}]\subseteq \mf{h}_{\overline{v}}$ but $\ell \not \in \widehat{\mf{k}}$. 
\end{prop}

\begin{remark}
   It is likely that the converse is also true: if $(z, \overline{v})\in F$ then there is a $\gamma (t)\subseteq G$ such that $y\stackrel{\gamma}{\rightarrow} (y, \overline{v})$. The proof is possibly through blow-ups. 
\end{remark}    
Let $\lambda (t)\subseteq G$ be a 1-PS and $V=\oplus_i V_i$ be the weight decomposition of $V$ with respect to $\lambda$. Let $V_{d+}=\oplus_{i\geq d} V_i$. Note that $V_{d+}$ is $P(\lambda)$-stable. 

\begin{lemma}
    \label{lemma:m[[t]]}
    Suppose that we have a $\gamma(t)\subseteq G$ such that $y\stackrel{\gamma}{\rightarrow}(z, \overline{y_1})$ of degree $1$, and $m(t)\in G({\cal O})$ is a 1-PS such that $m(0)=e_G$, then $y\stackrel{m\cdot \gamma}{\longrightarrow}(z, \overline{y_1})$ as well. 
\end{lemma}

\noindent
{\bf Proof}: Let $\gamma (t) y= y_0 +y_1 t^1+\ldots $, with $z=y_0$. Then $y'(t)=m(t)\gamma (t) y=m(t)y_0 +m(t)y_1 t^1 +\ldots$. Note that since $m(t)\in G({\cal O})$, the degree zero term $y'(0)$ is of course, $m(0)y_0$. Now coming to the tangent of approach, we see that $\mf{v}(m\cdot \gamma)=\mf{v}(m)+m(t) \mf{v}(\gamma) m(t)^{-1}$. Thus, we have the equation:
\[ (\mf{v}(m) + m(t)\mf{v}(\gamma) m(t)^{-1})y'(t)=\dot{y'}(t)=y'_1 +ty'_2+\ldots\]
But 
\[ m(t)\mf{v}(\gamma) m(t)^{-1}y'(t)=m(t)\mf{v}(\gamma)y(t)= m(t) \dot{y}(t)\]
Again, since $m(t)\in G({\cal O})$, $\mf{v}(m)\in \mf{g}\otimes {\cal O}$. The degree $0$ term in $\mf{v}(m)y'(t)$ is $\mf{v}(m)(0)y'(0)$. 
We thus have $y'_1=\mf{v}(m)(0)y_0 + m(0)y_1\in y_1 +T_z $. This proves the claim. $\Box$

\textcolor{red}{to be eaten.}
\begin{prop}    \label{prop:pathsimplify}
    Let $\lambda\subseteq T\subseteq G$. Suppose that we have a $\gamma(t)\subseteq G$ such that $y\stackrel{\gamma}{\rightarrow}(z, v)$, then there is an $p'_0\in \overline{P}(\lambda), \tau \in W(G;T)$ and a $p(t)\in P(\lambda)$ such that $y\stackrel{p}{\rightarrow}(\tau p'_0 z, \tau p'_0 v)$.
\end{prop}

\noindent{\bf Proof}: By the Iwasawa decomposition, we have $\gamma (t)=m[[t]]\cdot D(t) \cdot b(t)$, where $m[[t]]\in G[[t]]$ and $m(0)\in G$, $D\subseteq T(t)$, with diagonal entries $t^d$, $d\in \Z$ and $b(t)\in P(\lambda)$. We may write $m[[t]]=m'[[t]] m_0$ such that $m'[[t]]\in G[[t]]$ and $m'(0)=e_G$ and $m_0\in G$. Then, we have $y\stackrel{m_0 D(t) b(t)}{\longrightarrow}(z, v)$. We may now write $m_0 $ as $m_0 =(p'_0)^{-1} \tau^{-1} p_1$, with $p'_0 \in \overline{P}(\lambda), p_1\in P(\lambda)$. The result follows by applying $\tau p_0$ to the limit. The 1-PF  $p(t)$ is $p(t)=p_1 D(t) b(t)$. $\Box$
\textcolor{red}{upto here.}

\begin{lemma}
    \label{lemma:pathdegree}
    Let $\gamma (t)$ be of the form $\gamma (t)=\mu(t) h(t)$ where $h(t)\in G({\cal O})$ such that $h(0)\in G$, and $\mu(t)\subseteq T$ is a non-trivial diagonal 1-PS. Then
    \begin{enumerate}
    \item $\mf{v}(\mu)=t^{-1} \ell$, where $\ell\in \mf{g}$ is a non-zero diagonal element. 
    \item Let $\mf{g}=\oplus_j \mf{g}^{\ell}_j$ be the weight space decomposition of $\mf{g}$ with respect to $\ell$. Then $\mf{v}(\gamma)(t)=\sum_{i\geq a} \lie{g}_i t^i + t^{-1} \ell$, where $\lie{g}_i\in \oplus_{k\leq i} \mf{g}_k$.  
    \end{enumerate}
\end{lemma}

\noindent
{\bf Proof}: For (i), note that if $\mu(t)=diag((t^{d_i}))$ then $\dot{\mu}(t)(\mu(t))^{-1}= ((d_i t^{-1}))=t^{-1} \ell$ with $\ell=diag((d_i))$. For (ii), note that $\mf{v}(\gamma)=\mu(t) \mf{v}(h)\mu(t)^{-1}+ \mf{v}(\mu(t))$. Since $h\in G({\cal O})$, $\dot{h(t)}h(t)^{-1}\in \mf{g}\otimes {\cal O}=\sum_{i\geq o} \lie{h}_i t^i$, with $\lie{h}_i\in \mf{g}$. Whence $\mf{v}(\gamma)=(\sum_{i\geq 0} \sum_j \lie{h}_i^j t^{i+j})+ t^{-1} \ell$, with $\lie{h}_i^j \in \mf{g}_j$. Thus, we have:
\[ \mf{v}(\gamma)= \sum_d t^d (\sum_j \lie{h}_{d-j}^j)+t^{-1} \ell=
\sum_d t^d \lie{g}_d+t^{-1} \ell\]
with $\lie{g}_d=\sum_j \lie{h}_{d-j}^j$. Now since $\lie{h}_k=0$ for $k<0$, we must have $d-j\geq 0$, i.e., $j\leq d$. This proves (ii). $\Box$

\begin{lemma}
\label{lemma:basic-ell'}
    Let $\ell' =\ell +\lie{g}_{-1}$ as above. Then:
    \begin{enumerate}
        \item The element $\ell'$ is semisimple and $\mu'(t)=t^{\ell'}= \bar{u}(t)\mu(t)\bar{u}(t)^{-1}$, where $\bar{u}(t)\in \overline{U}(\mu)({\cal O})$.
        \item Let $ \mf{g}=\mf{g}^{\ell'}_i$ be the weight space decomposition of $\mf{g}$ with respect to $\ell'$. Then for every $i$ we have $\sum_{k\leq i} \mf{g}^\ell_k = \sum_{k\leq i} \mf{g}^{\ell'}_k$.
        \item Let $ V=V^{\ell'}_i$ be the weight space decomposition of $V$ with respect to $\ell'$. Then for every $i$ we have $\sum_{k\leq i} V^\ell_k = \sum_{k\leq i} V^{\ell'}_k$.
    \end{enumerate}
\end{lemma}

\noindent{\bf Proof}: To be written.

\begin{prop}
    \label{prop:hu}
    Let $y\stackrel{\gamma}{\rightarrow} z$ be such that $\gamma (t)= \mu(t) h(t)$ as above. Let $y(t)=y_0 t^0 +\dots $ and $\mf{v}(\gamma)= t^a \lie{g}_a + \ldots $ such that $a<0$. Then:
    \begin{enumerate}
        \item $\lie{g}_a \in \mf{h}$ and is nilpotent. Moreover, there exists a $\lie{u}\in \mf{h}$ such that $[\lie{g}_a , \lie{u}]=0$. 
        \item Let $u=exp(\lie{u})\in H$ and let $\gamma '(t)=u\gamma(t)$. Then
        \[ \mf{v}(\gamma')= \lie{g}'_bt^b +\ldots\]
        with $b>a$, and 
        \[ y'(t)= y_0 + tuy_1 + \ldots \]
        Thus $y\stackrel{\gamma'}{\rightarrow} (y_0 , uy_1 )$.
    \end{enumerate}
\end{prop}

\noindent
{\bf Proof}: Part (1) follows from earlier lemmas and general theory. For part (2), note that $\mf{v}(\gamma') = u\mf{v}(\gamma)u^{-1}$. Now the choice of $u$ implies that $u\lie{g}_a u^{-1}=0$. \textcolor{red}{Duh! It is itself. But others may get hiked...so perhaps some point.} $\Box$

\begin{defn}
    For a limit $y\stackrel{\gamma}{\rightarrow} (z,v)$, we say $v$ is a semisimple tangent of approach if $\widehat{\mf{v}(\gamma)}v=\alpha v$ with $\alpha \neq 0$.  
\end{defn}
\begin{lemma}
\label{lemma:ygree}
    Let $\gamma (t)$ be of the form $\gamma (t)=\mu(t) h(t)$ as above. Let $\mf{v}(t)=\mf{v}(\gamma)=\sum_{i\geq a} \lie{g}_i t^i +t^{-1} \ell$ be as in Lemma \ref{lemma:pathdegree}. Suppose that $\gamma(t)y=y_0 +ty_1 + \ldots +t^k y_k +\ldots $ is a semisimple limit, i.e., $\widehat{\mf{v}} \overline{y}_1 = \alpha \overline{y}_1$, with $ \overline{y}_1\neq 0$ and $\alpha \neq 0$. Then (i) $a= -1, \alpha=1$ and $\widehat{\mf{v}}=\ell+\lie{g}_{-1}$ where $\lie{g}_{-1}\in \oplus_{i\leq -1} \mf{g}_i$, (ii) let $\ell' = \ell +\lie{g}_{-1}$, then $\ell' y_0=0$ and $\ell' \overline{y_1}=\overline{y_1}=\overline{y^1_1}$, where $y^1_1$ is the $\ell'$-weight $1$ component of $y_1$. (iii) $y_k \in \oplus_{i\leq k} V^{\ell'}_i $ for all $k\geq 0$.
\end{lemma}

\noindent
{\bf Proof}: If $a<-1$ then $\widehat{\mf{v}}=\lie{g}_a$ is nilpotent, which would contradict the hypothesis that $y_1$ is a non-trivial eigenvector. Next $\lie{g}_{-1}\in \oplus_{k\leq -1} \mf{g}_k $ and thus, since $\ell$ is of weight zero, $\widehat{\mf{v}}=\ell +\lie{g}_{-1}$ is non-zero. This proves (i) except for the claim that $\alpha =1$, which follows from (ii) below. 
Coming to (ii), the equation $\mf{v}y(t)=\dot{y}(t)$ gives us:
\[ \begin{array}{rcl}
(\ell +\lie{g}_{-1})y_0 &=& 0 \\
(\ell +\lie{g}_{-1})y_1 + \lie{g}_0 y_0&=& y_1 \\
\vdots & & \vdots \\
(\ell +\lie{g}_{-1})y_i + \lie{g}_0 y_{i-1}+\ldots + \lie{g}_{i-1} y_0 &=& iy_i \\
\end{array}
\]
The second equation and the hypothesis $\widehat{\mf{v}} \overline{y}_1 = \alpha \overline{y}_1$ shows that $\alpha =1$, thus completing the proof of (i).
Let us now prove (ii). If $y_0 = \sum_{i\leq a} y_0^i$ is the $\ell '$-weight space decomposition of $y_0$, then we have $\ell' y_0 = \sum_{i\leq a} iy_0^i$. Now, separating the first equation by $\ell'$-weights gives us $a=0$ and $y^i_0 =0$ for $i<0$. Thus $y_0 = y^0_0$. Thus $\lie{g}_0 y_0 = \sum_{k\leq 0} \lie{g}_0^k y_0 $ is a weight space decomposition. Coming to the second equation, and assuming $y_1 =\sum_{i\leq b} y^i_1$, we have for all $i\leq b$:
\[ \begin{array}{rcl}
iy^i_1  &=& y_1^i \mbox{ for $i \geq 1$} \\ 
iy^i_1 + \lie{g}_0^i y_0 &=& y_1^i \mbox{ for $i\leq 0$} 
\end{array}
\]

This gives us that $b=1$ and that $y^i_1=0$ for $i>1$. For $i<1$, we have $((1-i)y_1^i = \lie{g}^i_0 y_0 $ is in the tangent space of $y_0$. Thus $\ell ' \overline{y_1}=\overline{y_1}$. This proves (ii). We prove (iii) by induction. 
Assuming the statement holds for $0,1,\ldots, k-1$, let us examine the equation:
\[ (\ell +\lie{g}_{-1})y_k + \lie{g}_0 y_{k-1}+\ldots + \lie{g}_{k-1} y_0 = ky_k\]
modulo $\oplus_{i\leq k-1} V_i$ to get:
\[ (\ell +\lie{g}_{-1}) y_k \equiv ky_k \: (mod \: \oplus_{i\leq k-1} V_i)\]
Thus, if $y_k =\sum_i y_k^i$ is the weight space decomposition of $y_k$, then we have $\ell' y_k =\sum_i i y_k^i$, and thus, we have that $y_k^i=0$ for $i>k$, which proves (iii) and the lemma. $\Box$

\begin{corollary}
\label{cor:firststep}
    Let $y\stackrel{\mu(t) h(t)}{\longrightarrow} (y_0 , y_1)$ be as above. Then there is a $y'\in O(y)$ and a 1-PS $\mu'(t)$ such that  $y\stackrel{\mu'(t)}{\longrightarrow} (y_0 , y_1)$.
\end{corollary}

\noindent
{\bf Proof}: Let us evaluate $\mu'(t^{-1})\gamma (t) y$. Since $y_i =\sum_{j\leq i} y^j_i$ is the $\ell'$-weight decomposition of $y_i$, we have:

\[ \begin{array}{rcl}
\mu'(t^{-1}) y(t) &=& \mu'(t^{-1}) (y_0 +y_1 t + \ldots  y_i t^i + \ldots )\\
&=& \sum_i \sum_j y_i^j t^{i-j}
\end{array}
\]
Note that $(i-j)\geq 0$ for all $i,j$ above and thus:
\[ y'=\lim_{t=0} \mu' (t)\gamma (t)y = y^0_0 +y^1_1  + \ldots +y^i_i \]
with $y_0^0 =y_0$ and $\overline{y_1^1}=\overline{y_1}$. 

Next note that 
\[ m(t)=\mu'(t^{-1})\gamma (t)=\bar{u}(t)\mu(t^{-1}) \bar{u}(t)^{-1} \mu(t) h(t) = \bar{u}(t)\bar{u}'(t) h(t)\in G({\cal O})\]
Thus $y' =m(0)y\in O(y)$. Finally 
\[ \mu'(t)y'= \sum_{i\geq 0} \mu'(t)y^i_i = y_0 +ty_1 +\ldots \]
And thus $y'\stackrel{\mu'}{\rightarrow} (y_0 , y_1)$ as before, with $\mu'$ as the required 1-PS. 
 $\Box$

\begin{prop}
    \label{prop:liedegree}
    Let $y\stackrel{\gamma}{\rightarrow} (y_0 , y_1)$ such that we have $\widehat{\mf{v}(t)} y_1 = \alpha y_1$ with $\alpha \neq 0$. Then there is an $m_0 \in G$ and a $\gamma'$ such that $y\stackrel{\gamma'}{\rightarrow}m_0 y'$ and a 1-PS $\mu (t)$ such that $m_0 y'\stackrel{\mu(t)}{\longrightarrow} (m_0 y_0 , m_0 y_1)$, with a semisimple $\ell \in G_{m_0 y_0}$ and $\ell m_0 y_1 = m_0 y_1$. Thus, upto conjugation by $G$,  $(y_0 , y_1)$ appears as a 1-PS limit of some point $y' \in \overline{O(y)}$. \textcolor{red}{Should this not be $y' \in O(y)$ itself?}
\end{prop}

\noindent
{\bf Proof}: Let us apply the Cartan-Iwahori decomposition to $\gamma (t)$ to get $\gamma (t)=m(t) \mu(t) u(t)$ with $\mu(t)\subseteq T$ as a 1-PS and $m(t), u(t)\in G[[t]]$. Let $m_0^{-1} =m(0)\in G$, then by Lemma \ref{lemma:m[[t]]}, we have $y\stackrel{\mu(t)u(t)}{\longrightarrow} (m_0 y_0 , m_0 y_1)$. Now, by Corollary \ref{cor:firststep}, there is an $ y'\in O(y)$ such that $y \stackrel{u(t)}{\longrightarrow} m_0 y'$ and  $y'\stackrel{\mu(t)}{\longrightarrow} (m_0 y_0 , m_0 y_1)$. This proves the proposition. $\Box$

\subsection{Stabilizer alignment}

\textcolor{blue}{This section investigates the question of alignment, i.e., the presence of a common semisimple stabilizer between $y$ and $z$.}
\noindent
\textcolor{red}{This lemma is lifted from JAMS-OnePS}

We begin with a technical lemma.

\begin{lemma} \label{lemma:uniorbit}
Let $y\stackrel{\lambda}{\rightarrow} z$ be as above and let $\ell$ be such that $t^{\ell}=\lambda(t)$. Then the following are equivalent:
\begin{enumerate}
    \item[(i)]$z\in O(y)$.
    \item[(ii)] There is an $\lie{k}\in {\cal K}$ such that $\widehat{\lie{k}}=\ell$.
    \item[(iii)] $y=uz$ for some $u\in U(\lambda)$.
    
\end{enumerate}
\end{lemma}

\noindent 
{\bf Proof}: Note that, $z$ is a degree $0$ limit of $y$ under $\lambda$. Therefore $\ell \in {\cal H}$ and $z$ is of $\ell$-weight $0$.

(i) $\implies$ (ii):  Since (i) holds, $dim({\cal H})=dim({\cal K})$ and ${\widehat{\cal K}} = {\cal H}$. Since $\ell \in {\cal H}$ and $\ell$ is homogeneous,  there must be a $\lie{k}$ such that $\widehat{\lie{k}}=\ell$. 

(ii) $\implies$ (iii):  Fix $\lie{k} \in {\cal K}$ such that $\widehat{\lie{k}}=\ell$. We can write
$\lie{k}=\ell +\lie{n}$
where $\lie{n} \in {\cal U}(\lambda)$.
Clearly, $\lie{k} \in {\cal P}(\lambda) \cap {\cal K}$.
The following claim can be easily established. We skip its proof.

\noindent
{\bf Claim}: $\lie{k} = \ell + \lie{n}$ is a semisimple element. Moreover, there is a $u\in U(\lambda)$ such that $u\lie{k}u^{-1}=\ell$. 

We complete the proof assuming the claim. With $y'=u.y$, we see that $y'\stackrel{\lambda}{\rightarrow}z$ and $\ell \in \mf{g}_{y'}$, the stabilizer of $y'$. Thus $\ell$ stabilizes both $y'$ and $z$. 

Now consider the $\ell$-weight decomposition of $y'$ as below
\[ y' =z_0 + z_1 +\ldots   \]
where $z_i \in V_i$ and $z_0 = z$ is the degree $0$ limit of $y$ under $\lambda$. 

Applying $\ell$ to the above equation, we get the following weight-space decomposition
\[ 0 = 0 + \ell z_1 + 2\ell z_2 + \ldots \]
This forces $z_i =0$ all $i \geq 1$ and $y'=uy=z$

(iii) $\implies $ (i) This is obvious. 
$\Box $

We now come to the main proposition.

\begin{prop}
    Let $z\in \overline{O(y)}$ and let the 1-PS $\lambda(t)$, stabilize $z$. Let $(-\lambda)$ denote that 1-PS $\lambda (t^{-1})$. Then either (A) there is a $y' \in (\overline{O(y)}-O(z))$ such that either (i) $\lambda$ stabilizes $y'$ or (ii) $y'\stackrel{\pm \lambda}{\longrightarrow} z$, or (B) there is a $\bar{\gamma }$. such that $\overline{\gamma} y= \sum_{+i\geq 0} t^i \overline{y}_i$ and each $\bar{y}_i$ is stabilized by $\lambda$. 
\end{prop}

\noindent{\bf Proof}: Since $z\in \overline{O(y)}$, there is a $\gamma$ such that $y\stackrel{\gamma}{\rightarrow} z$. Let $\gamma(t)y= \sum_{i\geq 0} t^i y_i$ and let $y_i=\sum_j y_i^j$ be the $\lambda$-weight decomposition of $y_i$. Note that by the hypothesis, $z=y_0 =y_0^0$. If $n=\max \{ |a| \mid V_i^{\lambda} \neq 0 \}$ and $N=n!$. 
Let $d_i = \max_j \{ y_i^j \neq 0 \}$ and $a_{\gamma}$ (or simply $a$) be $ \max_{i>0} \{ d_i/i \}$. Moreover, let $I=\{  i | i>0 \: and \: d_i/i=a \} \cup \{ 0\}$. Similarly, let $e_i = \min_j \{ y_i^j \neq 0 \}$ and $b_{\gamma}$ (or simply $b$) be $ \min_{i>0} \{ d_i/i \}$. Moreover, let $J=\{  i | i>0 \: and \: d_i/i=b \} \cup \{ 0\}$. If $a<0$ or $b>0$ then $y\stackrel{\pm\lambda}{\rightarrow} y_0=z$ and the proposition is proved. If $a=b=0$ then $y(t)=\gamma(t)y$ is stabilized by $\lambda$ and therefore (B) holds. That too proves the proposition. The cases which remain are when either $a>0$ or $b<0$ or both. 

Consider the case of $a>0$. Let $a=p/q$ with $\gcd(p,q)=1$ and $q=i$ and 
$\gamma' (t)=\lambda(t^{-q
})\gamma (t^p)$.  Then note that $I$ is a finite set and we have:
\[ \begin{array}{rcl}
\gamma' (t) y &=&  \sum_{i\geq 0} \sum_j t^{(pi-qj)} y_i^{j}\\
&=& t^0(\sum_{i\in I} y^{d_i}_i)+ \sum_{i>0} \sum_{j: pi-qj>0} t^{(pi-qj)} y_i^j
\end{array} \] 
Let $y'=\sum_{i\in I} y^{d_i}_i$. Then clearly $y' \in \overline{O(y)}$ and $y' \stackrel{\lambda}{\rightarrow} y_0$. If $y' \not \in O(z)$, then the proposition is proved. 

Consider the case when $y' \in O(z)$. Let $y'(t)=\sum_{i\in I} y_i^{d_i} t^i $. By Lemma \ref{lemma:uniorbit} below, we have $y'=u z$ for some $u\in U(\lambda)$. 

We use the following claim:

\noindent
{\bf Claim}: There is $u(t)\in U(\lambda)(\C[t])$ such that $u(t)\cdot z=y'(t)$. 

We outline a proof when $\lie{u}$ is of pure $\lambda$-weight.

Let $u=exp(\lie{u})$, where $\lie{u}\in \mf{u}(\lambda)$. Let $q=\min(I)$ and $p=d_q$, then the condition $y'=uz$ implies that $\lie{u}$  is of weight $p$. In other words, for any $v\in V_i^{\lambda}$, $\lie{u} v\in V_{i+p}^{\lambda}$.  

Let $u(t)=exp(t^q \lie{u})$. Then we have 
\[ u(t) y_0^0=y'(t)=\sum_{i\geq 0} t^{qi} y_{qi}^{pi} = \sum_{j\in I} t^j y_i^{d_j}\]
Note that for any $t^i v^j$, with  $v^j \in V_j^{\lambda}$, we have:
\[ u(t)^{-1} t^i v^j = t^i v^j +\sum_{r>0} t^{i+qr}v^{j+pr}\]
for some suitable $v^{j+qr}\in V_{j+qr}^{\lambda}$. Thus if $L_d = \{ (i',j')\in \Z^2| pi'-qj'=d\} $ and $\delta =pi-qj$ then all terms $t^{i'} v^{j'}$ in $u(t)^{-1}v^i t^j $ are such that $(i',j')\in L_{\delta}$. 

Let us examine $\gamma^1 (t) = u(t)^{-1}  \gamma (t)$. We have:
\[ \gamma^1 (t) y= t^0 y_0 + \sum_{d>0} \sum_{(i', j')\in L_d , i'>0} 't^{i'} y_{i'}'^{j'}\]
for some suitable $y_{i'}^{j'}\in V^{\lambda}_{j'}$

{\bf Statement ${\cal O}$}: We note that $y^1=\gamma^1 (t)y\in V\otimes {\cal O}$ and $a_{\gamma^1} < a_{\gamma}$ . \eat{Also note that $1/a$ is an integral power of $t^{1/N}$ and thus $\gamma^1 \in \C((t^{1/N}))$. Since $exp(t^{1/a} \lie{u})= \sum_{i\geq 0} \frac{t^{i/a}}{i!} \lie{u}^i$, we see that for all powers $t^j$ with $j<1/a$, the coefficients of $\gamma (t)$ and $\gamma^1 (t)$ agree. Finally, $\gamma^1 (t) y = t^0 y_0 + \sum_{j>0} (y^1)_j \cdot t^j$, where $j$ ranges over integers. } Moreover, the coefficients $y_j$ in $\gamma(t)y$ and $(y^1)_j$ in $\gamma^1 (t)y$ of $t^j$ agree for $j<q$. 

If we had chosen to tackle $b_{\gamma}=b<0$, then our $\gamma^1 (t)$ would be $u'(t)^{-1} \gamma (t)$, where $u(t)=exp(t^{q}\lie{u}')$, for a suitable $\lie{u}'\in \bar{\mf{u}}(\lambda)$, the Lie algebra of $\bar{U}(\lambda)$, the {\it opposite} unipotent group associated with $\lambda$. This would have given us $b_{\gamma^1} >b_{\gamma}$.

Thus, we may construct a sequence:
\[  \gamma \stackrel{u_1 (t)}{\rightarrow} \gamma^1 \stackrel{u_2 (t)}{\rightarrow} \ldots \stackrel{u_k (t)}{\rightarrow} \gamma^k \stackrel{u_{k+1} (t)}{\rightarrow} \ldots \]
where $u_i (t)\gamma^{i-1}(t)=\gamma^i (t)$. This sequence
either (i) terminates and we arrive at a suitable $y'$ proving the proposition, or (ii) by the argument in Statement ${\cal O}$, the sequence $(\gamma^k )\rightarrow \bar{\gamma}$ and so does the sequence $\gamma^k (t) y\rightarrow \bar{y}(t)$. We see then that if $\bar{y}(t)=\sum_{j\geq 0} \bar{y}_j t^j $, then $\bar{y}^0_0 =y_0$, and $\bar{y}_j^i=0$ for all $i\neq 0$. This implies that $\bar{y}(t)$ is stabilized by $\lambda$ and therefore (B) is satisfied. This completes the proof of the proposition. $\Box$

\begin{remark}
The final $\bar{\gamma}$ is the matrix product of $(\prod_{i\geq 1} \circ \: u_i (t) )\circ \gamma(t) $, where each $u_i (t)$ is a matrix with polynomial entries. We conjecture that if $\gamma (t)$ is analytic at $0$, then so is $\bar{\gamma}(t)$.    
\end{remark}

We now come to an important corollary on the boundary of the determinant.

\begin{corollary} \label{cor:codim1ps}
Let $V=Sym^n (X_n^*)$ under the action of $G=GL_{n^2} (\C)$ and $y=det_n \in V$. Let $\partial (y)=\overline{O(y)}-O(y)$. It is well know that $\partial (y)= X_1 \cup \ldots \cup X_r$ where each $X_i$ is an algebraic variety of co-dimension $1$. Suppose that $X=X_i=\overline{O(z)}$ for some point $z\in V$. Then there is a 1-PS $\lambda$ and a $y'\in O(y)$ such that $y'\stackrel{\pm \lambda}{\rightarrow}z$ or $\lambda (t)\in G_{y'}$. 
\end{corollary}

\noindent
{\bf Proof}: Let $G'=SL_{n^2} (\C)$. It is well known that $\overline{O_{G'}(z)}$ is itself a cone over $0$. Thus, there is a $\mu (t)\subseteq G'$ such that $\mu(t)z=t^d z$, for some $d\neq 0$. Consider $Z(t)=tI$, the center of $G$. Then, $Z(t)z=t^n z$ and thus there is a $\lambda \in G-G'$ such that $\lambda(t)z=z$ and thus $\lambda\in G_z$. 
Hence, by the above proposition, since there is no orbit $y'$ between $y$ and $z$, we must have $y'\in O(y)$. This proves the corollary. $\Box$

}

\bibliography{references}

@article{iwahori1965some,
  title={On some {B}ruhat decomposition and the structure of the {H}ecke rings of $ p $-adic {C}hevalley groups},
  author={Iwahori, Nagayoshi and Matsumoto, Hideya},
  journal={Publications Math{\'e}matiques de l'IH{\'E}S},
  volume={25},
  pages={5--48},
  year={1965}
}

@book{humphreys2012introduction,
  title={Introduction to Lie algebras and representation theory},
  author={Humphreys, James E},
  volume={9},
  year={2012},
  publisher={Springer Science \& Business Media}
}

@book{humphreys2012linear,
  title={Linear algebraic groups},
  author={Humphreys, James E},
  volume={21},
  year={2012},
  publisher={Springer Science \& Business Media}
}

@misc{ein2008jetschemessingularities,
      title={Jet schemes and singularities}, 
      author={Lawrence Ein and Mircea Mustata},
      year={2008},
      eprint={math/0612862},
      archivePrefix={arXiv},
      primaryClass={math.AG},
      url={https://arxiv.org/abs/math/0612862}, 
}

@article{ASS25jolt,
title ={Stabilizer Limits and Tangent Ideals},
author={Adsul, Bharat and Sohoni, Milind and Subrahmanyam, K V},
journal={Journal of Lie Theory (accepted for publication)},
year ={2026}
}

@misc{stacks-047N,
  author = {{Stacks Project Authors}},
  title = {{The Stacks Project}: Cartier's theorem: smoothness in characteristic zero},
  howpublished = {\url{https://stacks.math.columbia.edu/tag/047N}},
  note = {Lemma 39.8.2, Tag 047N}
}

@misc{stacks-02H6,
  author       = {{Stacks Project Authors}},
  title        = {{The Stacks Project}, Lemma 37.11.7:
                  Infinitesimal Lifting Criterion},
  howpublished = {Tag 02H6},
  url          = {https://stacks.math.columbia.edu/tag/02H6}
}

@book{mumford1994geometric,
  title={Geometric {I}nvariant {T}heory},
  author={Mumford, David and Fogarty, John and Kirwan, Frances},
  volume={34},
  year={1994},
  publisher={Springer Science \& Business Media}
}
\bibliographystyle{alpha}
\end{document}